\documentclass[preprint,12pt]{elsarticle}

\usepackage{lipsum}
\usepackage{amsfonts}
\usepackage{graphicx}
\usepackage[caption=false]{subfig}
\usepackage{enumerate}

\usepackage{epstopdf}
\usepackage{algorithmic}

\usepackage{url}
\usepackage{amsmath}
\usepackage{amssymb}
\usepackage{bm}

\makeatletter

\ifpdf
  \DeclareGraphicsExtensions{.eps,.pdf,.png,.jpg}
\else
  \DeclareGraphicsExtensions{.eps}
\fi

\newtheorem{cor}{Corollary}

\newcommand{\vf}{\bm{f}}
\newcommand{\vx}{\bm{x}}

\newcommand{\vtheta}{\bm{\theta}}

\newcommand{\vu}{\bm{u}}

\newcommand{\vY}{\bm{Y}}
\newcommand{\vw}{\bm{w}}
\newcommand{\vW}{\bm{W}}

\newcommand{\vb}{\bm{b}}
\newcommand{\vxi}{\bm{\xi}}

\newcommand{\vgamma}{\bm{\gamma}}
\newcommand{\vB}{\bm{B}}

\newcommand{\vk}{\bm{k}}
\newcommand{\vr}{\bm{r}}

\newcommand{\mW}{\bm{W}}
\newcommand{\mX}{\bm{X}}
\newcommand{\sR}{\mathbb{R}}

\newcommand{\fL}{\mathcal{L}}
\newcommand{\fN}{\mathcal{N}}
\newcommand{\fF}{\mathcal{F}}
\newcommand{\fS}{\mathcal{S}}

\newcommand{\ri}{\mathrm{i}}
\newcommand{\D}{\mathrm{d}}
\newcommand{\Exp}{\mathbb{E}}
\newcommand{\I}{\mathrm{i}}
\newcommand{\T}{\intercal}
\newcommand*\diff{\mathop{}\!\D}
\newcommand{\norm}[1]{\lVert#1\rVert}
\newcommand{\ReLU}{\mathrm{ReLU}}

\usepackage[markup=default,authormarkup=none]{changes}
\definechangesauthor[name={Reviewer 1}, color=purple]{R1} 
\definechangesauthor[name={Reviewer 2}, color=blue]{R2} 

\usepackage{cancel}

\journal{Journal of Computational Physics}

\begin{document}

\begin{frontmatter}



\title{On understanding and  overcoming spectral biases of deep neural network learning methods for solving PDEs} 

\author[university1]{Zhi-Qin John Xu \fnref{equal}} 
\ead{xuzhiqin@sjtu.edu.cn}

\author[university1]{Lulu Zhang \fnref{equal}} 
\ead{zhangl9661@sjtu.edu.cn}

\author[university2]{Wei Cai\corref{Corresponding}} 
\ead{cai@smu.edu}

\cortext[Corresponding]{Corresponding author.}
\fntext[equal]{Co-first author.}

\affiliation[university1]{organization={Institute of Natural Sciences and School of Mathematical Sciences, Shanghai Jiao Tong University},
            city={Shanghai},
            postcode={200240},
            country={China}}

\affiliation[university2]{organization={Dept of Mathematics, Southern Methodist University},
            city={Dallas},
            postcode={75252},
            state={Texas},
            country={USA}}
\begin{abstract}
In this review, we survey the latest approaches and techniques developed to overcome the spectral bias towards low frequency of deep neural network learning methods in learning multiple-frequency solutions of partial differential equations. Open problems and future research directions are also discussed.

\end{abstract}


\end{frontmatter}



\section{Introduction}



With its remarkable successes in image classification \cite{krizhevsky2012imagenet}, speech recognition \cite{hinton2012deep}, and natural language processing \cite{Lecun1998Gradient, bengio2009learning, LeCun2015deep}, the deep neural network (DNN) has recently also found many applications in scientific computing. These include modeling and predicting physical phenomena with various applications in material sciences, drug design, etc. \cite{weinan2020machinelearning,han2018deep, weinan2021algorithms}. Many of these applications involve solving partial differential equations (PDEs) such as the Schrodinger equation for quantum systems \cite{hermann2020deep, pfau2020ab}, Maxwell's equations for electromagnetic phenomena \cite{lim2022maxwellnet}, the Poisson-Boltzmann equation for molecular solvation \cite{liu2020multi,caibook}, the Navier-Stokes equations in fluid dynamics \cite{raissi2019physics, jin2021nsfnets}, the elasticity dynamic equations in civil engineering and geophysics \cite{Yin2022elasticity}, the Hamilton-Jacobi-Bellman equations from optimal controls \cite{weinan2017deep}, to list just a few. Traditionally, solutions to those equations in low dimensions are computed by finite element \cite{ciarlet2002fem,zienk2005fem,jiang1998least,bochev2009least}, finite difference \cite{strikwerda2004finite}, spectral \cite{hussaini1987spectral,Shen2011Spectral} or integral equation methods \cite{hack1995integral} with much success and maturity. However, these methods  face challenges when handling high-dimensional problems arising from treating stochastic effects, and material parameterization, or problems with complex geometries. Due to the compositional structure of the DNNs and the meshless random sampling during their learning, algorithms based on DNNs provide the potential of tackling high-dimensional problems without the curse of dimensionality \cite{han2018deep, weinan2021algorithms} where the traditional methods will not be able to handle at all,  or of not suffering the costs and difficulties of meshing complex geometries.


Many physical problems in scientific computing engender solutions with a wide range of frequencies. Examples include turbulent flows \cite{lumley72,chorin13} where energy interactions among different scales in the form of  vorticities and boundary layer effects; high-frequency waves in light-matter interactions from laser physics; and lithography and optical proximity correction in microchip designs at nanoscales and modern telecommunications \cite{saleh2019fundamentals, born2013principles,caibook}; quantum wave functions in many body electron systems \cite{landau2013quantum}; Wigner quantum transport in nano-electronics \cite{rammer2018quantum,caibook}; combustion dynamics \cite{zhang2022multi}, among many others. These wide-range frequency phenomena are described by various PDEs mentioned above.
The multigrid  (MG) methods \cite{brandt1977multi, xu1992iterative} developed to address the challenges of multiple-frequency issues in solving PDEs, though mainly in the low dimensional $R^3$ space, took advantage of the fast high-frequency smoothing property of classical iterative methods (e.g., Jacobi and Gauss-Seidel) based on the splitting of the coefficient matrices from the discretization of the PDEs. The matrices from finite element and finite difference methods are sparse, thus making the iterative solvers \cite{saad2003iterative} such as conjugate gradient and GMRES attractive, however, ill-conditioned with conditional numbers at the order of $\frac{1}{h^2}$ ( and even worse for spectral methods \cite{hussaini1987spectral}) for elliptic problems where the mesh size $h$ is required to be fine enough to resolve the highest frequency component of the solution based on Shannon sampling theory \cite{shannon1948mathematical}. Consequently, solving these linear systems necessitates the use of preconditioners \cite{ciaramella2022iterative,olver2013fast} with these iterative solvers.
The main reason is that the large disparity in the spectra of the coefficient matrices creates difficulties in removing errors across all frequencies using  traditional iterative solvers.  However, the multigrid method with V- and W-cycles of smoothing and corrections between different grids can produce an efficient solution in achieving uniform reduction of errors over all possible frequencies, and after many decades of research,  the MG methods now provides a frequency uniform convergence for solving elliptic type of PDEs.

On the other hand, DNN as a new computational model for solving PDEs is still in its early stage of development. It lacks the sophistication, efficiency, reliability, and robustness of traditional numerical methods, even with recent intensive research in this area. The research community faces many challenges, one of them, the focus of this review, is the well-known spectral bias (or frequency principle) nature of the DNNs  \cite{xu2019frequency,rahaman2019spectral,luo2022upper,xu2022overview}. This bias is evident in the approximation of functions or PDE solutions with a wide range of frequency content. Typically, the training of the DNNs  demonstrates a preference towards the low frequency before the content of higher frequencies can be learned at latter stage of the training.
This behavior contrasts with that of the classic iterative solvers. DNN-based computational methods developed for PDEs, particularly for high-dimensional problems common in many-body physics or problems with stochasticity in most engineering applications, or in complex geometry problems, face unique challenges due to the spectral bias.
Thus, in order to make DNN-based algorithms to reach the same robustness and practicality as the classic MG method even for high dimensional problems, it is imperative to overcome the shortcoming of the spectral bias in learning wide-range frequency solutions. Intensive researches have been carried out  by the computational community to improve the performance of DNNs.
This review represents a first effort to summarize the state-of-the-art research directed towards this goal and aims to inspire further research to resolve the issue of spectral bias of DNN learning methods for PDEs more effectively.

The following lists background materials and surveyed methods, developed to overcome spectral biases of DNN based learning methods for PDEs:

\begin{itemize}

\item Preliminaries on DNN learning of PDE solutions

\item Spectral bias of deep neural networks

\item Frequency shifting approach

\item Frequency scaling approach

\item Hybrid approach

\item  Multi-scale neural operator and diffusion models.


\end{itemize}

\section{Preliminaries on DNN learning of PDE solutions}
\subsection{Deep neural network (DNN)}
We introduce some conventional notations for DNNs.
A $L$-layer neural network is defined recursively as,

\begin{equation}
    \begin{aligned}
         & \vf_{\vtheta}^{[0]}(\vx)=\vx,                                                                                      \\
         & \vf_{\vtheta}^{[l]}(\vx)=\sigma\circ(\mW^{[l-1]} \vf_{\vtheta}^{[l-1]}(\vx) + \vb^{[l-1]}), \quad 1\leq l\leq L-1, \\
         & \vf_{\vtheta}(\vx) \equiv \vf_{\vtheta}^{[L]}(\vx)=\mW^{[L-1]} \vf_{\vtheta}^{[L-1]}(\vx) + \vb^{[L-1]},
    \end{aligned}
\end{equation}
where $\vx\in \sR^{d}$, $\mW^{[l]} \in \sR^{m_{l+1}\times m_{l}}$, $\vb^{[l]}=\sR^{m_{l+1}}$, $m_0=d_{\rm in}=d$ is the input dimension, $m_{L}=d_{\rm o}$ is the output dimension, $\sigma$ is a scalar nonlinear function, and ``$\circ$'' means entry-wise operation.
We denote the set of all parameters in the network by $\vtheta$.


The approximation capabilities of DNNs have been investigated extensively.
Cybenko \cite{cybenko1989approximation} demonstrated a universal approximation theorem of neural networks that any continuous function can be arbitrarily well approximated by sufficiently wide DNNs with one hidden layer with a continuous sigmoidal nonlinearity activation function. The sigmoidal function needs to saturate as the input variable goes to infinity. Almost simultaneously, Hornik et al. \cite{hornik1989multilayer} proved an universal approximation theorem for squashing activation functions, which needs to be non-decreasing and also saturate as the input variable goes to infinity. Furthermore, Hornik \cite{hornik1991approximation} extended the activation functions to those arbitrarily bounded and non-constant ones. More generally, Leshno et al. \cite{leshno1993multilayer} illustrated that DNNs with a locally bounded piecewise continuous activation function could approximate any continuous function to any degree of accuracy if and only if the network's activation function is nonlinear and not a polynomial. E et al. \cite{e2022Barron} showed that functions from a Barron space can be approximated by two-layer neural networks without the curse of dimensionality.  Lu et al. \cite{lu2021deep} established a quantitative approximation error characterization of deep ReLU networks for smooth functions in terms of arbitrary width and depth simultaneously. These collective works underscore the powerful approximation ability of DNNs. However, in practical applications, the use of appropriate optimization methods is also crucial for finding satisfactory solutions.

\subsection{DNN based methods for solving PDEs}

In solving PDEs, the empirical loss can be defined without labels as an unsupervised learning. In such cases, all possible data can be sampled with a non-zero probability, different from classical supervised learning problems.
Let us consider the following  equation,
\begin{equation}\label{pde}
        \begin{split}
        \mathcal{L}[u] (\vx) &=f(\vx), \quad \vx \in \Omega, \\
        \Gamma[u] (\vx)&=g(\vx), \quad \vx \in \partial \Omega, \\
        \end{split}
\end{equation}
where $\Omega \subset \mathbb{R}^d$ is a domain and $\partial \Omega $ is its boundary.
Note that time-dependent problems can also be formulated as equations of this form with the time variable included as an additional dimension (though,  only initial condition(s) at $t=0$ will be needed as an initial boundary value problem).
Also, when dealing with time-dependent problems, it is crucial to consider the causality of the system, which refers to the causal connections between consecutive physical events. In this context, one event, the cause, must precede another event, the effect, in time.
Taking causality into account when solving time-dependent problems is important, as it can significantly improve the performance and accuracy of the solutions \cite{liu2022causalitydeeponetcausalresponseslinear}.

\subsubsection{PDE residual based learning -- physics-informed neural network (PINN)}
A common way to build a loss function borrows the idea of minimization of residuals from spectral methods
\cite{hussaini1987spectral}
and the least square finite element (LSFE) method \cite{bochev2009least,jiang1998least}.
Along this line, the earlier work of neural network methods \cite{dissanayake1994neural,lagaris98,lagaris00,sirignano2018dgm,berg2018unified} and later the physics-informed neural network (PINN) \cite{raissi2019physics,kharazmi2021hp} and deep least square methods using first order system formulations of PDEs \cite{caiz2020,wang2020multi} attempt to find the solution of differential equations (DEs) by minimizing the least squares of the PDE's residuals.  And, the empirical
risk with a least square loss for the PDE \eqref{pde} can be defined by
\begin{equation}
\widetilde{L}_{pde}
(\vtheta)=
\frac{1}{n}\sum_{i=1}^{n}\Vert
\fL[u_{\vtheta}](\vx_{i}) - f(\vx_{i})
\Vert_{2}^{2},\label{pinn_loss}
\end{equation}
where the dataset $\{\vx_{i}\}_{i=1}^{n}$ is randomly sampled from $\Omega$ at
each iteration step and, in addition, a boundary loss will be used to enforce the boundary condition, i.e., for the
Dirichlet BC,
\begin{equation}
L_{bdry}(\vtheta)=\frac{1}{N_{bdry}}%
{\displaystyle\sum\limits_{k=1}^{N_{bdry}}}
|u_{\vtheta}(\vx_{k})-g(\vx_{k})|^{2},\quad \vx_{k}\in\partial \Omega.
\label{bdryloss}
\end{equation}
where the dataset $\{\vx_{k}\}_{k=1}^{N_{bdry}}$ is also randomly sampled from $\partial \Omega$ at
each iteration step.
Note that at each iteration, adaptive sampling based on $\widetilde{L}_{pde}(\vtheta)$ and $L_{bdry}(\vtheta)$ might lead to more accurate results \cite{yan2023failureinformed, yan2024failureinformed}, however, in this paper, for simplicity of presentation, a uniform distribution is chosen for the sampling.

\subsubsection{Variational energy based Deep-Ritz method}
Another very natural candidate for the loss function $L(\vtheta)$ for the learning optimization procedure is the Ritz energy from the variational principle of elliptic operators~\cite{weinan2018deep}.

The Deep-Ritz method \cite{weinan2018deep} seeks a variational solution $u(\vx)$  through the following minimization problem,
\begin{equation}
u=\arg\min_{v\in H^1_0(\Omega) }J(v),\label{Ritz}%
\end{equation}
where $H^1_0$ is the Sobolev space for the set of admissible functions (also called trial functions, here represented by $v$).

A typical variational problem is based on the following functional
\begin{equation}
    J(v)=\int_{\Omega}(\frac{1}{2}|\nabla v(x)|^2 - f(x)v(x))dx,
    \label{varloss}
\end{equation}
where $f$ is a given function, representing external forcing to the system under consideration.

The PINN requires high-order derivatives in the definition of the loss function, which could potentially accelerate the convergence of high-frequency components during training as differentiation in the frequency domain is equivalent to a multiplication with the wave number $k$ in the frequency domain, therefore increasing the weight of high-frequency components in the loss function. Lu et al. \cite{lu2019deepxde} and E et al. \cite{ma2020machine} have provided examples showing that differentiation accelerates the convergence of high frequencies. However, higher-order differentiation could also lead to training instability.   In contrast, the DeepRitz, which reduces the differentiation by one order for second-order elliptic PDEs, offers a more stable training, albeit at the cost of slower convergence for high-frequency components.
Generally, these methods still exhibit the characteristics of a preferential convergence in low-frequency components,  and
these formulations have been used for multi-frequency problems.
 The current literature lacks a rigorous and systemic comparison of the pros and cons of various formulations for different types of problems.

\section{Spectral bias of deep neural networks}
 In the learning process of DNNs, there exists
a spectral bias (studied in \cite{xu_training_2018,xu2019frequency} as a frequency principle) during training, supported by much empirical evidence \cite{xu_training_2018,xu2019frequency,rahaman2019spectral,xu2022overview} and many theoretical works \cite{xu2019frequency,luo2019theory,luo2022exact,zhang2021linear,cao2019towards,bordelon2020spectrum,luo2022upper}. The spectral bias refers to the fact that  DNNs tend to fit target functions from low to high frequencies, contrary to the behavior of conventional iterative  schemes of numerical linear algebra \cite{Axelsson}. As a result of the spectral bias, DNNs have been shown to face significant challenges in learning solutions with a wide range of frequencies.

\subsection{Empirical evidence of spectral bias}

To gain insight into the training process, we visualize the training process in the frequency domain with a  one-dimensional data with some given frequencies. The training samples are extracted from a one-dimensional target function with three frequencies $f^{*}(x)=\sin(x)+\sin(3x)+\sin(5x)$ over the interval $[-3.14, 3.14]$ in an evenly spaced range.  The samples are represented as $\{x_{i},f(x_{i})\}_{i=0}^{n-1}$. We compute the discrete Fourier transform (DFT) of $f^{*}(x)$ and the DNN output $f_{\vtheta}(x)$, and denote them as $\hat{f}^{*}{(k)}$ and $\hat{f}_{\vtheta}{(k)}$, respectively, where $k$ denotes the frequency. Fig.~\ref{fig:onelayer}(a) illustrates the target function with its three distinct frequencies (indicated by red dots in  Fig.~\ref{fig:onelayer}(b)). The DNN is trained by using a mean squared error (MSE) loss and we consider the relative error of  each frequency, $\Delta_{F}(k)=||\hat{h}_{k}-\hat{f}_{k}||/||\hat{f}_{k}||$.
Fig.~\ref{fig:onelayer}(c) shows the relative errors for the process, indicating that the DNN learns in the order of low to high frequencies.

The spectral bias of DNNs has been widely observed in high-dimensional problems and diverse settings in a series of experiments conducted in \cite{xu2019frequency}. These experiments used common optimization methods, different network structures, various loss functions, and more. In addition, the spectral bias has also been found in the training of non-gradient based algorithms \cite{ma2021frequency}, such L-BFGS and quasi-Newton's methods,  etc.

\begin{center}
\begin{figure}
\begin{centering}
\subfloat[$f^{*}(x)$]{\includegraphics[width=0.3\textwidth]{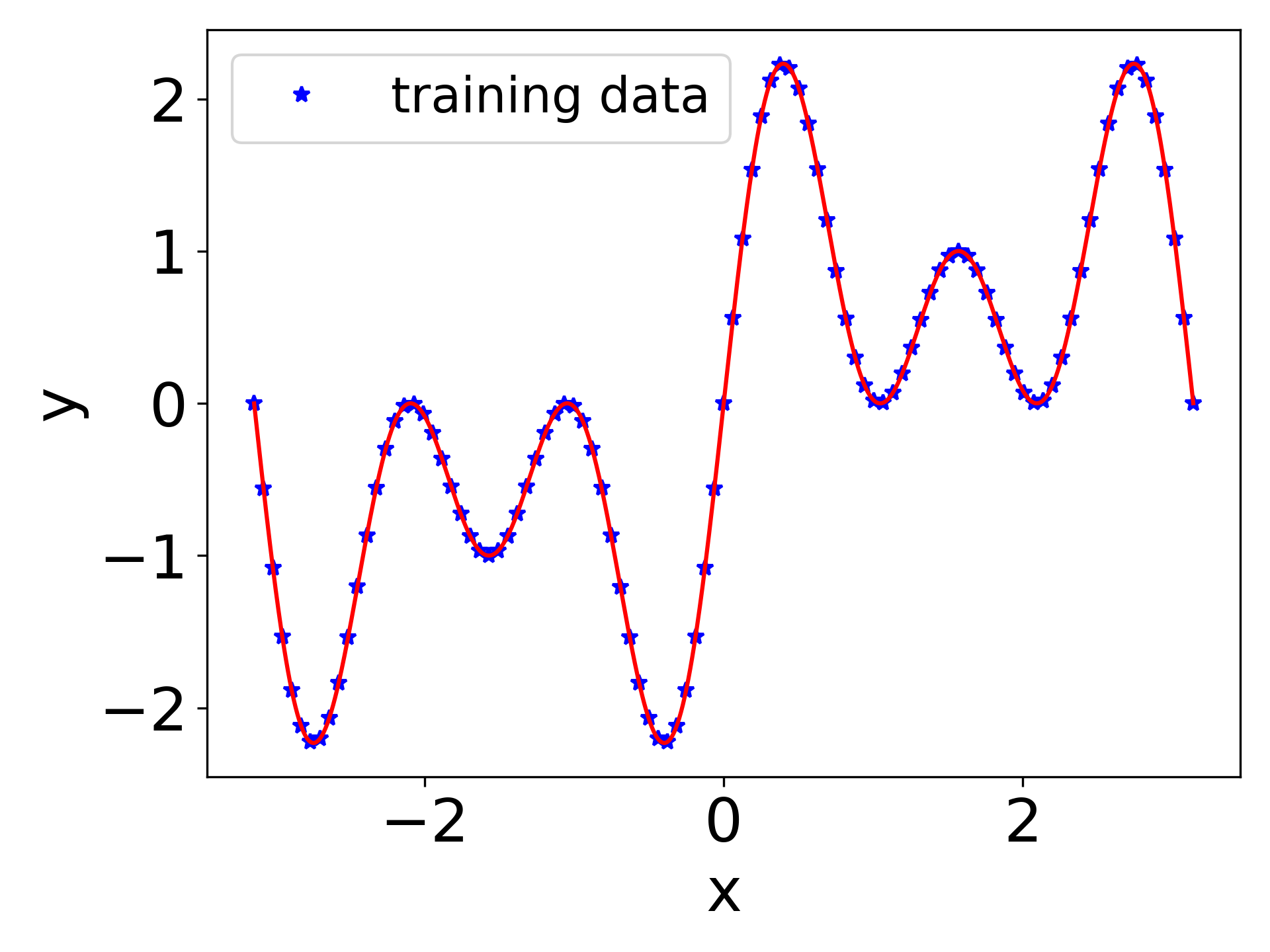} }
\subfloat[$||\hat{f}^{*}(k)||$]{\includegraphics[width=0.3\textwidth]{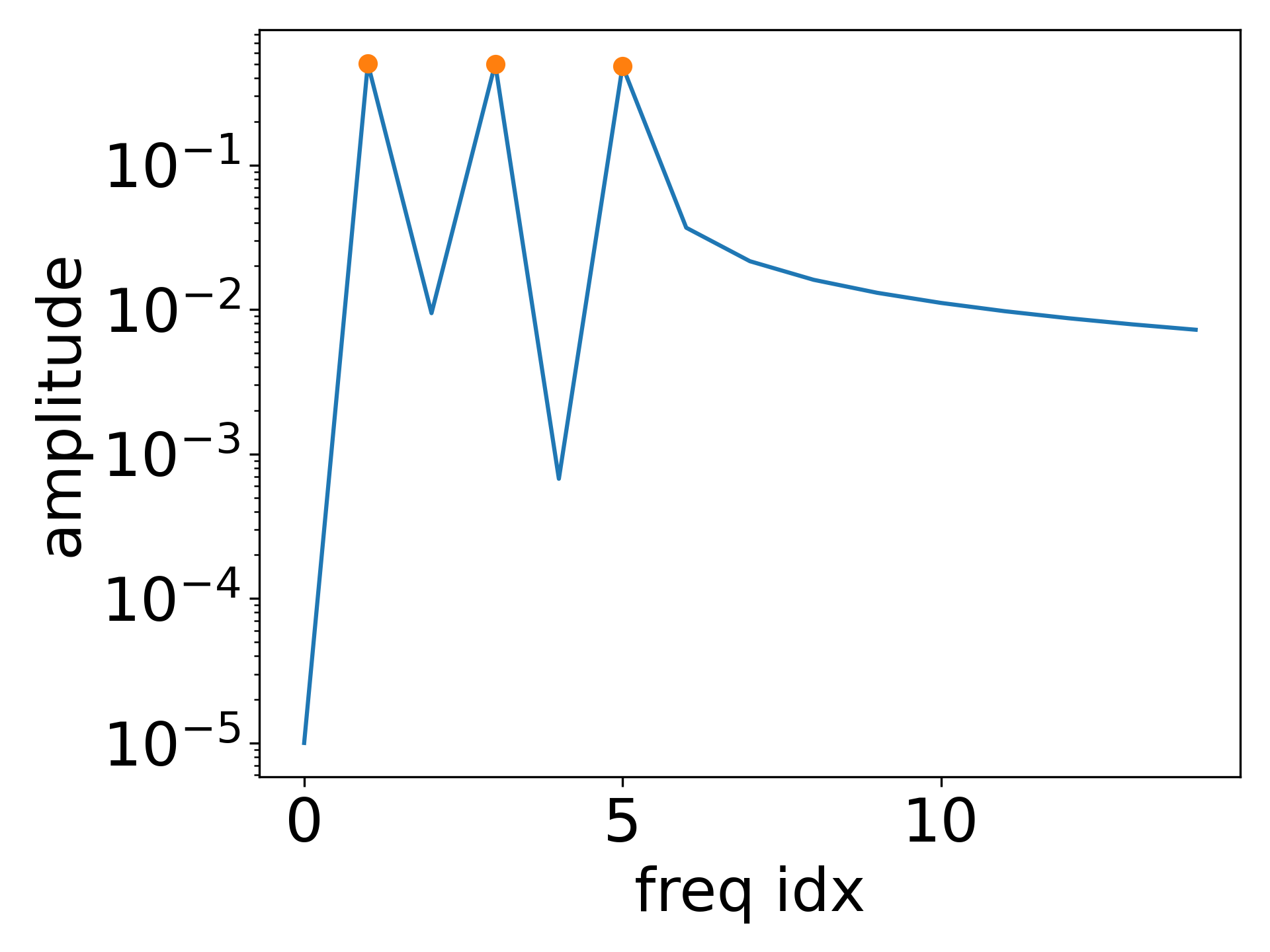} }
\subfloat[relative error]{\includegraphics[width=0.3\textwidth]{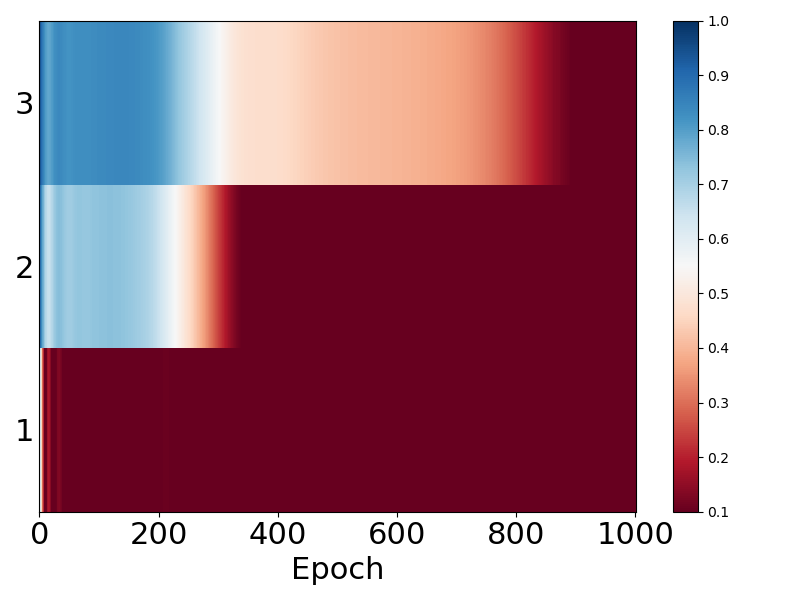} }

\par\end{centering}
\caption{1-d
frequency principle.
(a) $f^{*}(x)$. (b) : $||\hat{f}^{*}(k)||$. (c) $\Delta_{F}(k)$
of three important frequencies (indicated by red dots in the inset
of (b)) against different training epochs (horizontal axis). Here the y-axis represents the 3 frequencies, and the colorbar indicates the loss values.
\label{fig:onelayer} }
\end{figure}
\par\end{center}

\subsection{Qualitative theory of spectral bias}
To gain an intuitive understanding of the underlying mechanism of the spectral bias, we review a simple analysis \cite{xu2018understanding,xu2019frequency,luo2019theory}. For this case, we consider the activation function
\begin{equation*}
    \sigma(x)= tanh(x)=\frac{e^{x}-e^{-x}}{e^{x}+e^{-x}},\quad x\in \mathbb{R},
\end{equation*}
and a DNN with one hidden layer, having $m$ neurons, a 1-dimensional input $x$, and a 1-dimensional output, i.e.,
\begin{equation}
    h(x)=\sum_{j=1}^{m}a_{j}\sigma(w_{j}x+b_{j}),\quad a_{j},w_{j},b_{j}\in{\rm \mathbb{R}},\label{eq: DNNmath}
\end{equation}
where $w_{j}$, $a_{j}$, and $b_{j}$ are training parameters. We also denote $\vtheta=\{\theta_{lj}\}$ with $\theta_{1j}=a_{j}$,
$\theta_{2j}=w_{j}$, and $\theta_{3j}=b_{j}$, $j=1,\cdots,m$.

Then, the Fourier transform of $h(x)$ can be computed as follows
\begin{equation}
    \hat{h}(k)=\sum_{j=1}^{m}\frac{2\pi a_{j}\I}{|w_{j}|}\exp \Big(\frac{\I b_{j}k}{w_{j}}\Big)\frac{1}{\exp(-\frac{\pi k}{2w_{j}})-\exp(\frac{\pi k}{2w_{j}})},\label{eq:FTW}
\end{equation}
where $k$ denotes the frequency. The loss at a frequency $k$ is $L(k)=\frac{1}{2}||\hat{h}(k)-\hat{f}^{*}(k)||^{2}$, where ${f}^{*}$ is the target function and the total loss function is defined as: $L=\int_{-\infty}^{+\infty}L(k)\diff{k}$. According to Parseval's theorem, this loss function in the Fourier domain is exactly the common mean squared error loss, that is, $L=\int_{-\infty}^{+\infty}\frac{1}{2}(h(x)-f(x))^{2}\diff{x}$.
Therefore, to analyse the error decay from a gradient descent training, we can use the loss in the frequency domain to see the error dependence on frequency. In practice, the weights are usually much smaller than $1$, then, the last term in Eq. (\ref{eq:FTW}) is approximately $\exp\left(-|\pi k/2w_{j}|\right)$. Therefore, the absolute contribution from frequency $k$ to the gradient for $\theta_{lj}$ is
\begin{equation}
    \norm{\frac{\partial L(k)}{\partial\theta_{lj}}}\approx ||\hat{h}(k)-\hat{f}^{*}(k)||\exp\left(-|\pi k/2w_{j}|\right)
    F_{lj} ,\label{eq:DL2}
\end{equation}
where $\vtheta_{j}\triangleq\{w_{j},b_{j},a_{j}\}$, $\theta_{lj}\in\vtheta_{j}$, $F_{lj}$ is an $O(1)$ function depending on $\vtheta_{j}$ and $k$.
The term $\exp\left(-|\pi k/2w_{j}|\right)$ shows that low frequency plays a dominant role for small training weights $w_{j}$, which leads to faster convergence of low-frequency components as manifested by the spectral bias of the DNNs.

The preceding analysis of the spectral bias relies on the form of the activation function, and most activation functions, such as tanh and ReLU, exhibit a similar decaying behavior in the frequency domain and, consequently, the spectral bias will also be readily observable. Also from this analysis, modifying the loss function to impose greater weight on high-frequency components can impact frequency convergence. One simple approach to mitigate spectral biases is incorporating derivatives of the DNN output with respect to the input in the loss function as the Fourier transform of $\nabla_{\vx} f_{\vtheta}(x_i)$ equals the product of the transform of $f_{\vtheta}(\vx_i)$ and frequency $k$, effectively prioritizing higher frequencies in the loss function, or contributing a factor $k$ in Eq. (\ref{eq:DL2}).
Approaches along this line can be found in \cite{xu2022overview,lu2019deepxde}. And, to accelerate the learning of high-frequency in fluid simulation,  the loss of gradient information was also used in \cite{biland2019frequency}.

\subsection{Linear quantitative theory of spectral bias} \label{sec:ntkdynamics}
A series of works analyzed the spectral bias for two-layer wide DNNs with the Neural Tangent Kernel (NTK) approach \cite{jacot2018neural,luo2021phase} with specific  sample distributions \cite{cao2019towards,basri2019convergence,bordelon2020spectrum} or any finite samples \cite{zhang2019explicitizing,zhang2021linear,luo2022exact}. Additionally,  in \cite{weinan2020machine} the spectral bias was studied  from the perspective of integral equation.

In this subsection, we review analysis work of the frequency dependence of the training error using the NTK approach together with linearization of network, i.e., a linear frequency principle (LFP) analysis for the spectral bias \cite{luo2022exact}.  And, we consider the following gradient-descent flow dynamics of the empirical risk $L_{S}$ of a network function $f_{\vtheta}(\cdot)=f(\cdot,\vtheta)$ parameterized by $\vtheta$ on a set of training data $\{(\vx_i,y_i)\}_{i=1}^{n}$
\begin{equation}
	\left\{
	\begin{array}{l}
		\dot{\vtheta}=-\nabla_{\vtheta}L_{S}(\vtheta),  \\
		\vtheta(0)=\vtheta_0,
	\end{array}
	\right.
\end{equation}
where
\begin{equation}
	L_{S}(\vtheta)
	= \frac{1}{2}\sum_{i=1}^n(f(\vx_i,\vtheta)-y_i)^2.
\end{equation}
Then, the training dynamics of the output function $f(\cdot,\vtheta)$ are
\begin{align}
	\frac{\D}{\D t}f(\vx,\vtheta) &= \nabla_{\vtheta}f(\vx,\vtheta)\cdot\dot{\vtheta}
	= -\nabla_{\vtheta}f(\vx,\vtheta)\cdot\nabla_{\vtheta}L_{S}(\vtheta) \nonumber \\
	&= -\nabla_{\vtheta}f(\vx,\vtheta)\cdot\sum_{i=1}^n \nabla_{\vtheta}f(\vx_i,\vtheta)(f(\vx_i,\vtheta)-y_i) \nonumber \\
	&= -\sum_{i=1}^n K_m(\vx,\vx_i)(f(\vx_i,\vtheta)-y_i),
 \label{errorEq}
\end{align}
where the NTK $K_m(\vx,\vx')(t)$ \cite{jacot2018neural} at time $t$ and evaluated at $(\vx,\vx')\in\Omega\times\Omega$ is defined by
\begin{equation}
	K_m(\vx,\vx')(t)=\nabla_{\vtheta}f(\vx,\vtheta(t))\cdot\nabla_{\vtheta}f(\vx',\vtheta(t)).
 \label{ntk}
\end{equation}

Now, consider a two-layer DNN
\begin{align}
	f(\vx,\vtheta)
	&= \frac{1}{\sqrt{m}}\sum_{j=1}^{m}a_{j}\sigma(\vw_{j}^{\T}\vx+b_{j})\label{eq: 2layer-nn},
\end{align}
where the vector of all parameters  $\vtheta$ denotes all network parameters,
initialized by a standard Gaussian distribution.
As $m\to \infty$, the following linearization around initialization
\begin{equation}
	f^{{\rm lin}}\left(\vx;\vtheta(t)\right)=f\left(\vx;\vtheta(0)\right)+\nabla_{\vtheta}f\left(\vx;\vtheta(0)\right)\left(\vtheta(t)-\vtheta(0)\right) \label{eq:linearNN}
\end{equation}
is an effective approximation of $f\left(\vx;\vtheta(t)\right)$, i.e.,  $f^{{\rm lin}}\left(\vx;\vtheta(t)\right)\approx f\left(\vx;\vtheta(t)\right)$ for any $t$ \cite{jacot2018neural,lee_wide_2019}, which defines the linear regime.  Note that,  $f^{{\rm lin}}\left(\vx;\vtheta(t)\right)$, linear in $\vtheta$ and nonlinear in $\vx$, reserves the universal approximation power of $f\left(\vx;\vtheta(t)\right)$ as $m\to \infty$. In the rest of this sub-section, we do not distinguish $f\left(\vx;\vtheta(t)\right)$ from $f^{{\rm lin}}\left(\vx;\vtheta(t)\right)$.

In the linear regime, one can prove  \cite{jacot2018neural}
\begin{equation}
    K^{*}(\vx,\vx'):=\lim_{m\rightarrow\infty} K_m(\vx,\vx')(t)=K(\vx,\vx')(0)
\end{equation}
for any $t$.
The gradient descent of the model thus becomes
\begin{equation}
	\frac{\D}{\D t}\vu(\vx,t)=-\sum_{i=1}^n K^{*}(\vx,\vx_i)\vu(\vx_i,t),
\end{equation}
where the residual $\vu(\vx,t)=f(\vx,\vtheta(t))-f^{*}(\vx)$, $f^{*}(\vx$ is the target function and $y_i = f^{*}(\vx_i)$. Denote $\mX\in \sR ^{n\times d}$ and $\vY\in \sR ^{n}$ as the training data, $ u(\mX):=u(\mX,\vtheta(t))\in \sR^{n},\nabla_{\vtheta}u(\mX,\vtheta(t))\in \sR^{n\times M}$ (M is the number of parameters), let $K^{*}=[K^{*}(\vx_i,\vx_j)]_{ij}\in\sR^{M\times M}$ with a slight abuse of the notation $K^{*}$, then,
\begin{equation}
	\frac{\D u(\mX)}{\D t}=- K^{*} u(\mX). \label{eq:eigenK}
\end{equation}

In a continuous form, one can define the empirical density $\rho(\vx)=\sum_{i=1}^n\delta(\vx-\vx_i)/n$ and further denote $u_{\rho}(\vx)=u(\vx)\rho(\vx)$. Therefore, the dynamics for $u$ becomes
\begin{equation}
	\frac{\D}{\D t}u(\vx,t)=-\int_{\sR^d}K^{*}(\vx,\vx')u_{\rho}(\vx',t)
	\diff{\vx'}.\label{eq..DynamicsFiniteWidth}
\end{equation}
This continuous form gives the integral equation analyzed in \cite{weinan2020machine}. The convergence analysis of the dynamics in Eq. (\ref{eq:eigenK}) can be done by performing eigen decomposition of $K^{*}$.

For the two-layer neural network (\ref{eq: 2layer-nn}), we assume that $b\sim\fN(0,\sigma_{b}^{2})$ with $\sigma_{b}\gg 1$. As a special case, we have the exact LFP dynamics for the cases where the activation function is ReLU.

\begin{cor}[Luo et al., (2022) \cite{luo2022exact}: LFP operator for ReLU] \label{cor..LFP.operator.for.ReLU}
    Under mild assumptions and if $\sigma_b\gg 1$ and $\sigma=\ReLU$, then the dynamics \eqref{eq..DynamicsFiniteWidth} satisfies the following expression,
    \begin{equation} \label{thmdyna.ReLU}
        \langle\partial_t\fF[u], \phi\rangle = -\left\langle \fL[\fF[u_{\rho}]], \phi \right\rangle+O(\sigma_b^{-3}),
    \end{equation}
    where $\phi\in \fS(\sR^d)$ is a test function and the LFP operator reads as
    \begin{equation} \label{eq..lfpoperatorthm.ReLU}
        \begin{aligned}
            \fL[\fF[u_{\rho}]]
             & = \frac{\Gamma(d/2)}{2\sqrt{2}\pi^{(d+1)/2}\sigma_b}\Exp_{a,r}\left[\frac{r^3}{16\pi^4\norm{\vxi}^{d+3}} + \frac{a^2 r}{4\pi^2\norm{\vxi}^{d+1}}\right]\fF[u_{\rho}](\vxi)                                     \\
             & \quad -\frac{\Gamma(d/2)}{2\sqrt{2}\pi^{(d+1)/2}\sigma_b}\nabla\cdot\left (\Exp_{a,r}\left[\frac{a^2 r}{4\pi^2\norm{\vxi}^{d+1}}\right]\nabla\fF[u_{\rho}](\vxi) \right),
        \end{aligned}
    \end{equation}
    where the expectations are taken w.r.t. initial parameter distribution, $r = \norm{\vw}$, $\fF[\cdot]$ indicates Fourier transform.
\end{cor}
These results explicitly estimate the convergence of each frequency, which confirms the observed spectral bias.

\section{Frequency shifting approach}

To address the spectral bias of DNN,  one approach is to uses frequency domain manipulation. This strategy
involves converting the higher frequency components of the data to lower frequencies prior to training, speeding up the learning process, and subsequently transforming the learned representations back to their original high-frequency range.

\begin{figure}[ptbh]
\centering
\includegraphics[width=0.5\linewidth]{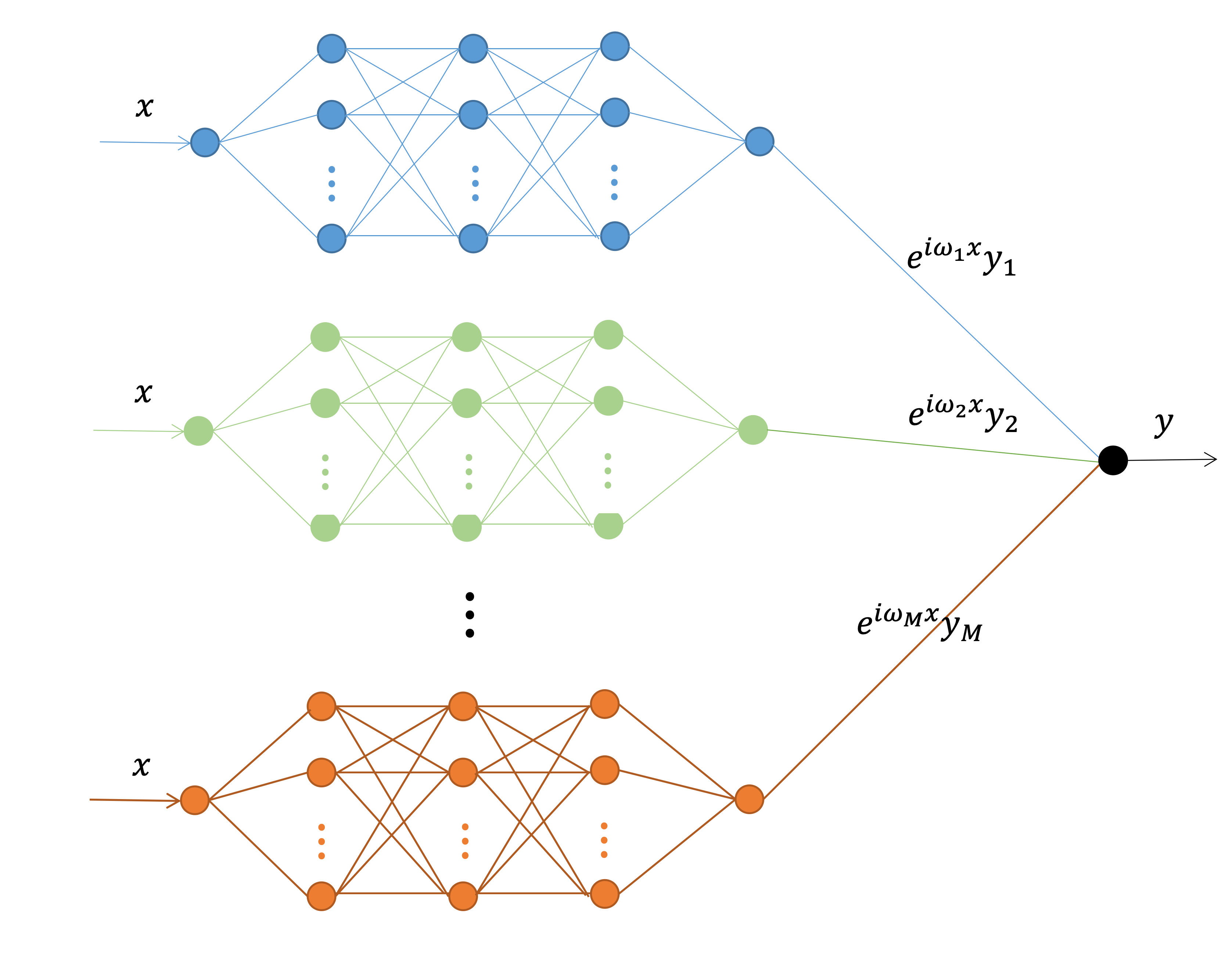}
\caption{
Illustration
of a PhaseDNN structure.}%
\label{PhaseDNN_schematic}%
\end{figure}

An example of this approach is the phase shift DNN (PhaseDNN) \cite{cai2019phasednn}.
The PhaseDNN was developed to tackle the challenges of high-frequency learning for regressions and solving PDEs such as Helmholtz equations. It transforms high-frequency data components into a lower-frequency spectrum by using phase shifts in the frequency domain. This enables the network to learn from these modified data quickly and, upon completion, map them back to the original high frequencies by reverse phase shifts. This method has been shown to be effective in achieving an uniform convergence across a broad frequency spectrum for wave propagation problems. A typical PhaseDNN architecture is shown in Fig. \ref{PhaseDNN_schematic}, which comprises a series of smaller sub-DNNs, each tasked with approximating a specific range of high-frequency content. With the help of phase shifts in the frequency domain, these sub-DNNs are capable of learning high-frequency information at a convergence rate typical of low-frequency problems, thus improving the overall performance and accuracy of the network in handling wideband high-frequency data.

To illustrate how the PhaseDNN arises, let us consider an one-dimensional problem for a band-limited function $f(x)$ within a frequency range
$[-K_{0},K_{0}]$, i.e.
\[
\text{Supp}\hat{f}(k)\subset\lbrack-K_{0},K_{0}]=[-M\Delta k,M\Delta k],\Delta
k=K_{0}/M.
\]

We first partition $[-K_{0},K_{0}]$  with uniform subintervals,
\begin{equation}
\lbrack-K_{0},K_{0}]\subset %
{\displaystyle\bigcup\limits_{j=-M}^{M}}
A_{j},A_{j}=[\omega_{j} -\frac{\Delta k}{2},\omega_{j} +
\frac{\Delta k}{2}],\omega_{j}=j\Delta k,j=-M,\cdots,M.
\end{equation}

We can decompose the target function $f(x)$ in the Fourier space as
follows,%
\begin{equation}
\hat{f}(k)=
{\displaystyle\sum\limits_{j=-M}^{M}}
\chi_{A_{j}}(k)\hat{f}(k)\triangleq%
{\displaystyle\sum\limits_{j=-M}^{M}}
\hat{f_{j}}(k),\label{k-decomp}%
\end{equation}
where $\chi_{A_{j}}$ is the indicator function for the interval $A_{j},$ which
will give a corresponding decomposition in $x$-space as%
\begin{equation}
f(x)=
{\displaystyle\sum\limits_{j=-M}^{M}}
f_{j}(x),\label{x-decomp}%
\end{equation}
where
\begin{equation}
    f_{j}(x)=\mathcal{F}^{-1}[\hat{f_{j}}](x).
\end{equation}

The decomposition (\ref{x-decomp}) involves $2M+1$ functions $f_{j}(x)$, whose
frequency spectrum is limited to $[\omega_{j} -\frac{\Delta k}{2},\omega_{j} +
\frac{\Delta k}{2}]$. Therefore, a
simple phase shift by $\omega_{j}$ could translate its spectrum to $[-\Delta k/2,\Delta k/2]$, allowing it to
be learned quickly by a relatively small DNN $T_{j}(x)$ with a few
training epoches. Namely, $\hat{f_{j}}(k+\omega_{j})$ is supported in
$[-\Delta k/2,\Delta k/2]$,  then its inverse Fourier transform $\mathcal{F}
^{-1}\left[  \hat{f_{j}}(k+\omega_{j})\right]$, denoted as
\begin{equation}
f_{j}^{\text{shift}}(x)=\mathcal{F}^{-1}\left[  \hat{f_{j}}(k+\omega
_{j})\right]  (x)\label{rjshift}%
\end{equation}
can be learned quickly by a small DNN $T_{j}(x,\vtheta_j)$ by minimizing a loss
function
\begin{equation}
L_{j}(\vtheta_j)=\int_{-\infty}^{\infty}|f_{j}^{\text{shift}}(x)-T_{j}%
(x,\vtheta_j)|^{2}d\vx\label{eq:Lj}%
\end{equation}
in some $n_{0}$-epoches of training.

Meanwhile, in the physical domain, we have
\begin{equation}
f_{j}^{\text{shift}}(x)=e^{-i\omega_{j}x}f_{j}%
(x).\label{dataShift}%
\end{equation}

Equation
(\ref{dataShift})
shows that once $f_{j}^{\text{shift}}(x)$ is learned,
$f_{j}(x)$ is also learned by simply removing the phase factor, i.e.,
\begin{equation}
f_{j}(x)\approx e^{i\omega_{j}x}T_{j}(x,\vtheta^{(n_{0})}).
\end{equation}
Now with all $f_{j}(x)$ for $-M\leq j\leq M$ learned after some $n_{0}$ steps of
training each, we have an approximation to $f(x)$ over the whole frequency range
$[-M\Delta k,M\Delta k]=[-K_0,K_0]$ as follows%
\begin{equation}
f(x)\approx T_{\vtheta}(x)= {\sum\limits_{j=-M}^{M}}e^{i\omega_{j}x}T_{j}%
(x,\vtheta_j),\label{PhaseDNN}%
\end{equation}
where $\vtheta=\{ \vtheta_j\}_{j=-M}^M$ is the parameters of the combined DNN, as depicted in Fig. \ref{PhaseDNN_schematic}.

The discussion above suggests that a
PhaseDNN
can be sought in the
following form (after converting to real functions),
\begin{equation}
T_{\vtheta}(x)=\sum_{m=0}^{M}A_{m}(x)\cos(\omega_{m}x)+B_{m}(x)\sin
(\omega_{m}x),
\end{equation}
where $A_{m}$,$B_{m}$ are DNNs with their individual network parameters and the shift frequencies $\omega_m$ in fact can be made as trainable parameters to best fit for the target functions.

\begin{figure*}[htpb]
\centering
\subfloat[]{\includegraphics[width=0.45\linewidth,trim=100 250 100 200, clip]{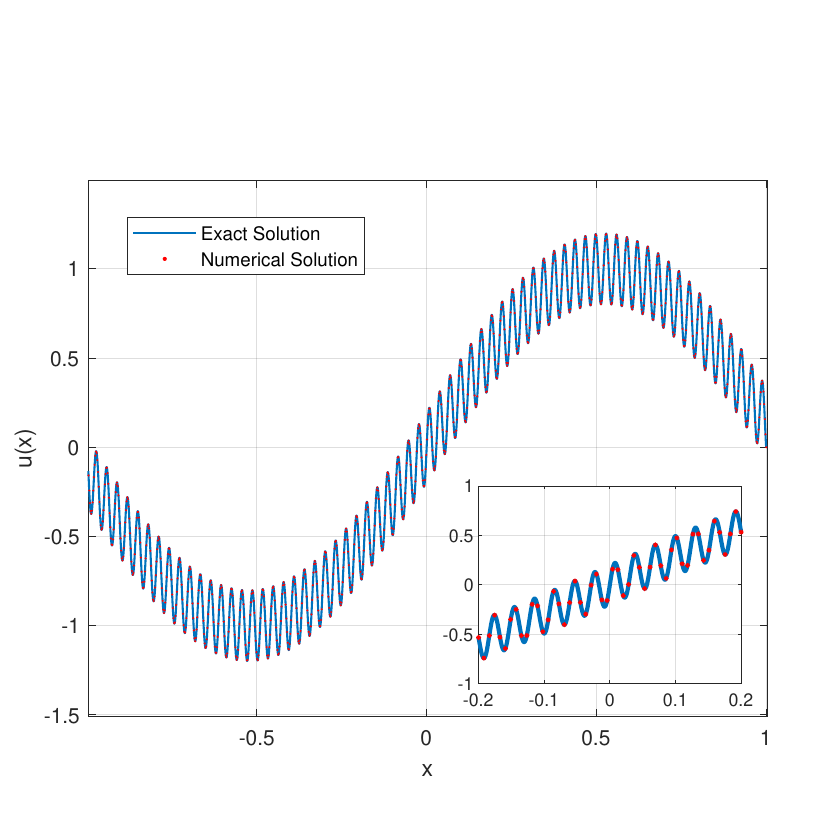}}
\subfloat[]{\includegraphics[width=0.45\linewidth,trim=100 250 100 200, clip]{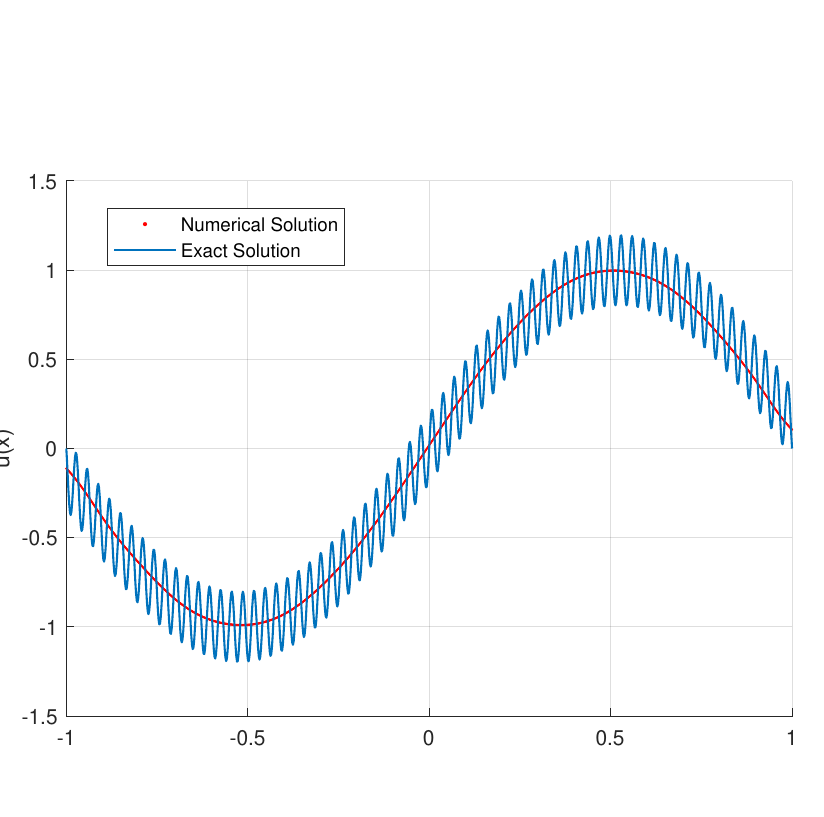}}
\caption{Results demonstrating the efficacy of using PhaseDNN  compared to a normal fully-connected DNN.
(a) Numerical and exact solution of Helmholtz equation using the PhaseDNN.
(b) Failure of a normal fully-connected DNN. }
\label{fig:add_PhaseDNN}
\end{figure*}

As an illustration, we use the PhaseDNN to solve the Helmholtz equation with a constant wave number,
\begin{equation}
u''+\lambda ^2 u = \sin \mu x,
\end{equation}
coupled with boundary conditions $u(-1)=u(1)=0$.
 For the PhaseDNN, the frequencies  $\{ \omega_m \}$ are selected to be 0, $\lambda$, and $\mu$.
Each $A_{m},B_{m}$ is set to be a 1-40-40-40-40-1 DNN.
The result for the case of $\lambda=205, \mu=3$ in Fig. \ref{fig:add_PhaseDNN} shows (right panel) that a normal fully-connected neural network failed to capture the high-frequency oscillation  while the PhaseDNN (left panel) gives an accurate approximation to the high-frequency solution. In practice, the frequency for the shift could be chosen as a hyperparameter for the training.

In \cite{peng2020accelerating}, a dictionary of $(1,\cos x,\sin x,\cos2x,\sin2x,\cdots,\cos kx,\sin kx)$
was used to construct a Prior Dictionary based Physics-Informed Neural Networks (PD-PINNs) to the same effect of  the PhaseDNN.
Also, to avoid using multiple neural networks for $A_{m}$'s and $B_{m}$'s,
a network with multiple outputs  was used to represent $A_{m}$'s
and $B_{m}$'s for the PhaseDNN \cite{kim2023phase}. Numerical experiments in
\cite{cai2019phasednn,peng2020accelerating,kim2023phase} demonstrate that the
PhaseDNN can be very effective in learning low-dimensional high-frequency
function. However, as the number of Fourier terms increases with the dimension
exponentially., the PhaseDNN suffers from the curse of dimensionality, and will face
difficulties for high-dimensional problems.

\section{Frequency scaling approach}

Another approach developed to overcome the spectral bias of DNN is to utilize some type of scaling in the frequency domain so that the higher frequency content will be converted to a lower frequency content before training. One method employing this approach is the multi-scale DNN (MscaleDNN) \cite{liu2020multi,zhang2023correction} using radial scaling in the frequency domain (thus equivalently in physical domain). Another related approach is the Fourier feature network \cite{tancik2020fourier}.  The activation function optimization method and the random feature method   \cite{chen2022bridging}  have also employed physical scaling to achieve multiscale effects.

\subsection{Multiscale DNN (MscaleDNN)}

\label{sec:mscalednn} The phase shift DNN (PhaseDNN) previously discussed will suffer from the curse of dimensionality for
high-dimensional problems due to the amount of phase shifts required along all coordinate directions. As another approach of frequency manipulation, the multiscale deep neural network (MscaleDNN)
proposed in \cite{caixu19,liu2020multi,zhang2023correction}  utilizes a radial scaling in the frequency domain, instead, to transform the frequency characteristics of the data,  as well as activation functions with localized frequency content such as functions with compact support or sin and cos functions. The radial scaling converts the problem of approximating high-frequency contents of target functions or PDEs' solutions to a problem of
learning about lower frequency functions, and compact support activation
functions  or simple sine or cosine functions facilitate the separation of frequency contents of the target
function to be approximated by individual DNNs.


\subsubsection{The idea of Multiscale DNN}

To illustrate how the MscaleDNN arises, again we consider a band-limited function $f(\mathbf{x})$, $\mathbf{x}
\in\mathbb{R}^{d}$, whose Fourier transform $\widehat{f}(\mathbf{k})$ has a
compact support, i.e.,%
\begin{equation}
\text{Supp}\widehat{f}(\mathbf{k})\subset B(K_{\text{max}})=\{\mathbf{k\in
}\mathbb{R}^{d},|\mathbf{k|\leq}K_{\text{max}}\}.\label{support}%
\end{equation}

We now partition the domain $B(K_{\text{max}})$ as union of $M$ concentric
annulus with uniform or non-uniform width, e.g., for the case of uniform
$K_{0}$-width
\begin{equation}
A_{i}=\{\mathbf{k\in}\mathbb{R}^{d},(i-1)K_{0}\leq|\mathbf{k|\leq}%
iK_{0}\},\quad K_{0}=K_{\text{max}}/M,\quad1\leq i\leq M\label{anulnus}%
\end{equation}
so that
\begin{equation}
B(K_{\text{max}})=%
{\displaystyle\bigcup\limits_{i=1}^{M}}
A_{i}.\label{part}%
\end{equation}

Now, we can decompose the function $\widehat{f}(\mathbf{k})$ as before
\begin{equation}
\widehat{f}(\mathbf{k})={\sum\limits_{i=1}^{M}}\chi_{A_{i}}(\mathbf{k}%
)\widehat{f}(\mathbf{k})\triangleq%
{\displaystyle\sum\limits_{i=1}^{M}}
\widehat{f}_{i}(\mathbf{k}),\label{POU}%
\end{equation}
where $\chi_{A_{i}}$ is the indicator function of the set $A_{i}$ and
\begin{equation}
\text{Supp}\widehat{f}_{i}(\mathbf{k})\subset A_{i}.\label{SupAi}%
\end{equation}

\bigskip The decomposition in the Fourier space give a corresponding one in
the physical space%
\begin{equation}
f(\mathbf{x})=
{\displaystyle\sum\limits_{i=1}^{M}}
f_{i}(\mathbf{x}),\label{Partx}%
\end{equation}
where
\begin{equation}
f_{i}(\mathbf{x})=\mathcal{F}^{-1}[\widehat{f}_{i}(\mathbf{k})](\mathbf{x}).\label{convol}%
\end{equation}

From (\ref{SupAi}), we can apply a simple down-scaling to convert the high
frequency region $A_{i}$ to a low-frequency one. Namely, we define a scaled
version of $\widehat{f}_{i}(\mathbf{k})$ as
\begin{equation}
\widehat{f}_{i}^{(\text{scale})}(\mathbf{k})=\widehat{f}_{i}(\alpha
_{i}\mathbf{k}),\qquad\alpha_{i}>1,\label{fkscale}%
\end{equation}
and, correspondingly in the physical space
\begin{equation}
f_{i}^{(\text{scale})}(\mathbf{x})=\frac{1}{\alpha_{i}^{d}}f_{i}(\frac{1}{\alpha_{i}%
}\mathbf{x}),\label{fscale}%
\end{equation}
or
\begin{equation}
f_{i}(\mathbf{x})=\alpha_{i}^{d}f_{i}^{(\text{scale})}(\alpha_{i}\mathbf{x}).\label{fscale_inv}%
\end{equation}
We can see the low-frequency spectrum of the scaled function $\widehat{f}%
_{i}^{(\text{scale})}(\mathbf{k})$ if $\alpha_{i}$ is chosen large enough,
i.e.,
\begin{equation}
\text{Supp}\widehat{f}_{i}^{(\text{scale})}(\mathbf{k})\subset\{\mathbf{k\in
}\mathbb{R}^{d},\frac{(i-1)K_{0}}{\alpha_{i}}\leq|\mathbf{k|\leq}\frac{iK_{0}%
}{\alpha_{i}}\}.\label{fssup}%
\end{equation}

Based on the low-frequency bias of DNNs, with $iK_{0}/\alpha_{i}$ being small, we can train a
DNN $f_{\vtheta^{n_{i}}}(\mathbf{x})$, with $\vtheta^{n_{i}}$ denoting the DNN
parameters, to learn $f_{i}^{(\text{scale})}(\mathbf{x})$ quickly
\begin{equation}
f_{i}^{(\text{scale})}(\mathbf{x})\sim f_{\vtheta^{n_{i}}}(\mathbf{x}%
),\label{DNNi}%
\end{equation}
giving an approximation to $f_{i}(\mathbf{x})$ immediately%
\begin{equation}
f_{i}(\mathbf{x})\sim\alpha_{i}^{d}f_{\vtheta^{n_{i}}}(\alpha_{i}%
\mathbf{x})\label{fi_app}%
\end{equation}
and to $f(\mathbf{x})$ as well
\begin{equation}
f(\mathbf{x})\sim%
{\displaystyle\sum\limits_{i=1}^{M}}
\alpha_{i}^{d}f_{\vtheta^{n_{i}}}(\alpha_{i}\mathbf{x}),\label{f_app}%
\end{equation}
which suggests the format of a multiscale DNNs (MscaleDNN) \cite{liu2020multi,zhang2023correction}. And the scale factors $\alpha_i$ can be made as trainable parameters to best fit for the target functions.
Two specific implementations of the MscaleDNNs are given below.

\begin{figure}[ptbh]
\centering
\subfloat[MscaleDNN-1]{\includegraphics[width=0.45\linewidth]{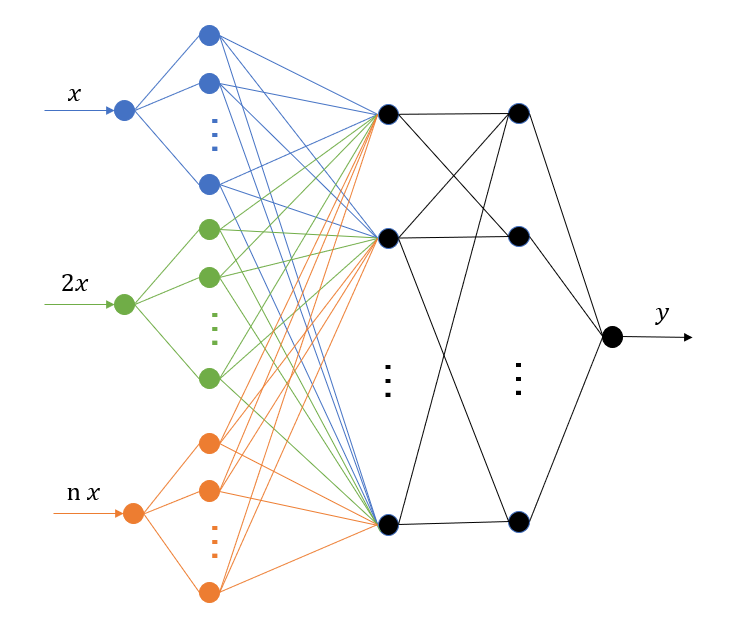}}
\subfloat[MscaleDNN-2]{\includegraphics[width=0.45\linewidth]{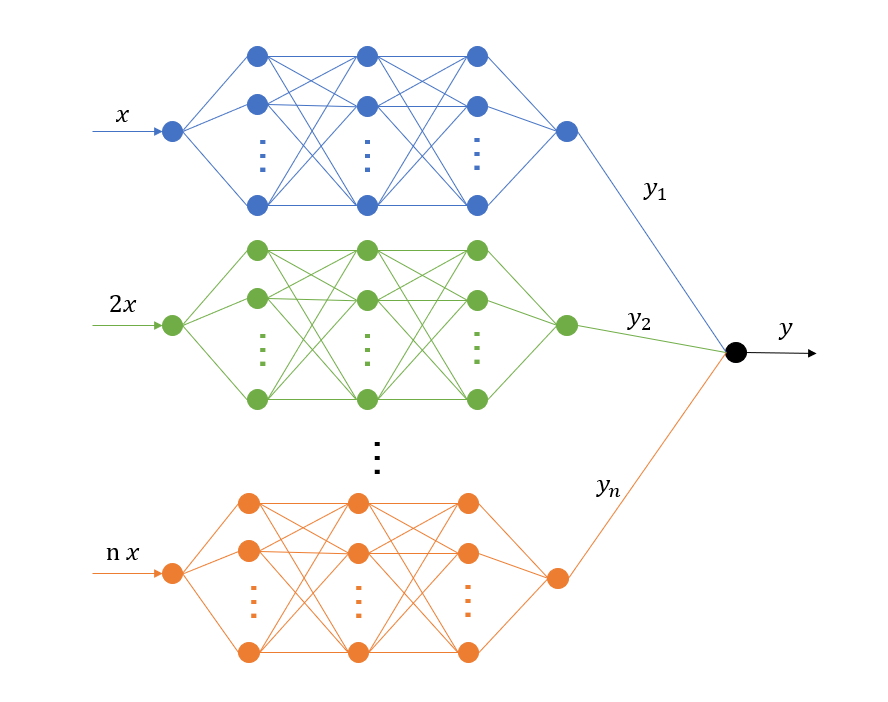}}\caption{Illustration
of two MscaleDNN structures \cite{liu2020multi}.}%
\label{net}%
\end{figure}


\medskip

\noindent\textbf{MscaleDNN-1} For the first kind, the neuron in the first
hidden-layer is grouped into to $N$ parts. The neuron in the $i$-th part
receives input $a_{i}x$, that is, its output is $\sigma(a_{i}\mathbf{W}%
\cdot\mathbf{x}+\mathbf{b})$, where $\mathbf{W}$, $\mathbf{x}$, $\mathbf{b}$
are weight, input, and bias parameters, respectively. A MscaleDNN takes the
following form
\begin{equation}
f_{\vtheta}(\mathbf{x})=\mathbf{W}^{[L-1]}\sigma\circ(\cdots(\mathbf{W}%
^{[1]}\sigma\circ(\mathbf{K}\odot(\mathbf{W}^{[0]}\mathbf{x})+\mathbf{b}%
^{[0]})+\mathbf{b}^{[1]})\cdots)+\mathbf{b}^{[L-1]},\label{mscalednn}%
\end{equation}
where $\mathbf{x}\in\mathbb{R}^{d}$, $\mathbf{W}^{[l]}\in\mathbb{R}%
^{m_{l+1}\times m_{l}}$, $m_{l}$ is the neuron number of $l$-th hidden layer,
$m_{0}=d$, $\mathbf{b}^{[l]}\in\mathbb{R}^{m_{l+1}}$, $\sigma$ is a scalar
function and \textquotedblleft$\circ$\textquotedblright\ means entry-wise
operation, $\odot$ is the Hadamard product and
\begin{equation}
\mathbf{K}=(\underbrace{a_{1},a_{1},\cdots,a_{1}}_{\text{1st part}}%
,a_{2},\cdots,a_{i-1},\underbrace{a_{i},a_{i},\cdots,a_{i}}_{\text{ith part}%
},\cdots,\underbrace{a_{N},a_{N}\cdots,a_{N}}_{\text{Nth part}})^{T},
\end{equation}
where $\mathbf{K}\in\mathbb{R}^{m_{1}}$, and the scale factor $a_{i}=i$ or $a_{i}=2^{i-1}$. This
structure is called Multi-scale DNN-1 (MscaleDNN-1), as depicted in Fig.
\ref{net}(a). \bigskip

\noindent\textbf{MscaleDNN-2}  A second kind of
multi-scale DNN is given in Fig. \ref{net}(b), as a sum of $N$ subnetworks, in
which each scaled input goes through a separate subnetwork. In MscaleDNN-2, weight
matrices from $\mathbf{W}^{[1]}$ to $\mathbf{W}^{[L-1]}$ are block diagonal.
Again, the scale factor $a_{i}=i$ or $a_{i}=2^{i-1}$.

Liu et al \cite{liu2020multi} has shown that the MscaleDNN can work well for high frequency fitting and PDE solving with a very complex geometry and the numerical experiments also found that the activation function in the first hidden
layer is important for the performance of MscaleDNNs, for example,
$\mathrm{sReLU}(x)=(x)_{+}(1-x)_{+}$ and $\phi(x)=(x-0)_{+}^{2}-3(x-1)_{+}%
^{2}+3(x-2)_{+}^{2}-(x-3)_{+}^{2}$ were used, where $x_{+}=\max\{x,0\}=\mathrm{ReLU}%
(x)$.
MscaleDNN with various activation functions is also explored, such as the one with sinusoidal activation functions used for PDEs with a point source \cite{huang2021solving}, the one with
residual connection and activation $\sigma(x)=0.1\ast(sinx+cosx)$ \cite{chen2023adaptive}, and the one with the ellipse Gaussian RBFs as activation function  \cite{Wang2023SolvingME}. In \cite{li2023subspace}, Li et al.  introduced an orthogonality
condition between low-frequency and high-frequency output subspaces through
a penalty term in the loss function, which is designed to
reduce the influences between low and high frequencies in the network, and achieved better
error in solving PDEs.

\begin{figure*}[htpb]
\centering
\subfloat[]{\includegraphics[width=0.45\linewidth]{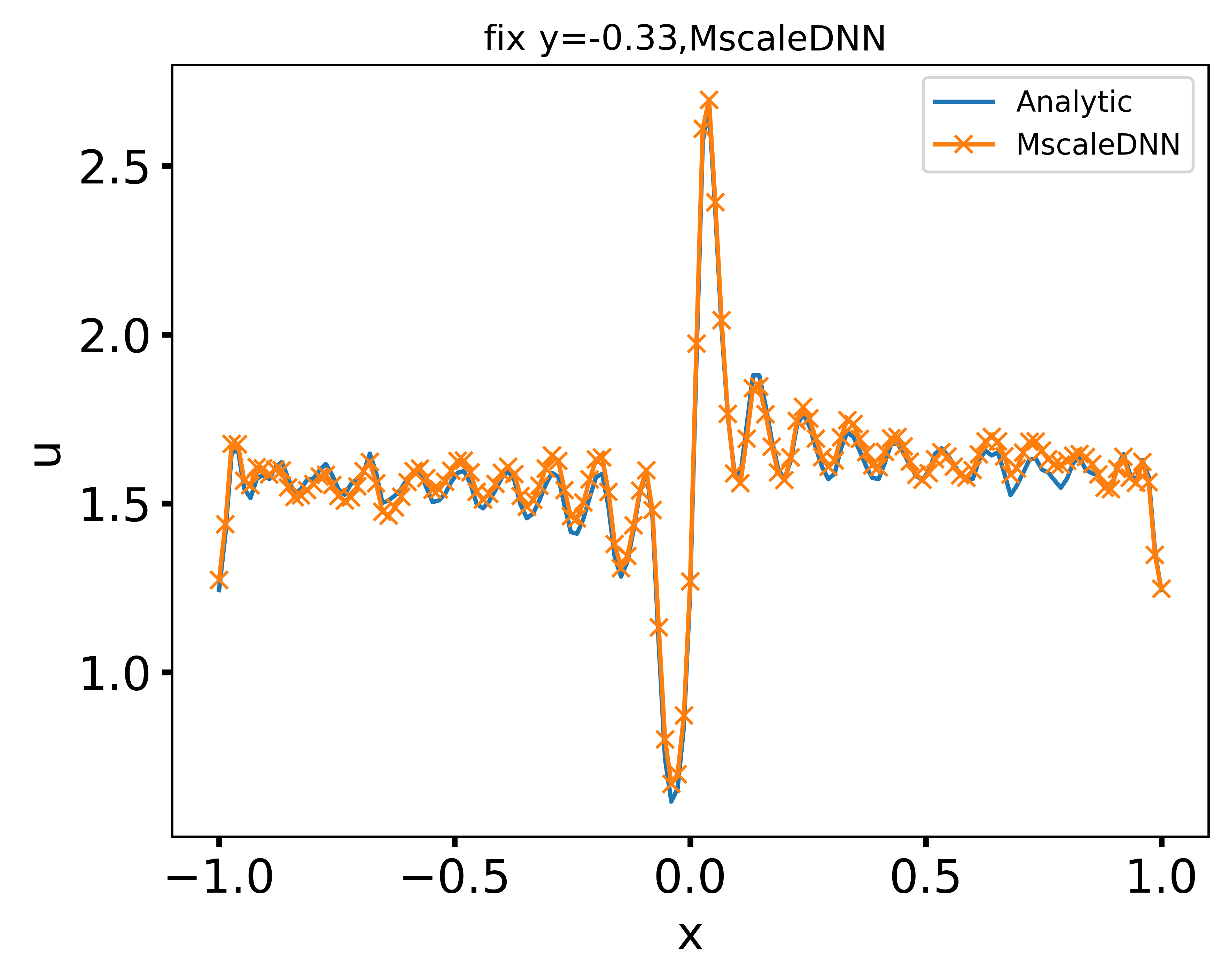}}
\subfloat[]{\includegraphics[width=0.45\linewidth]{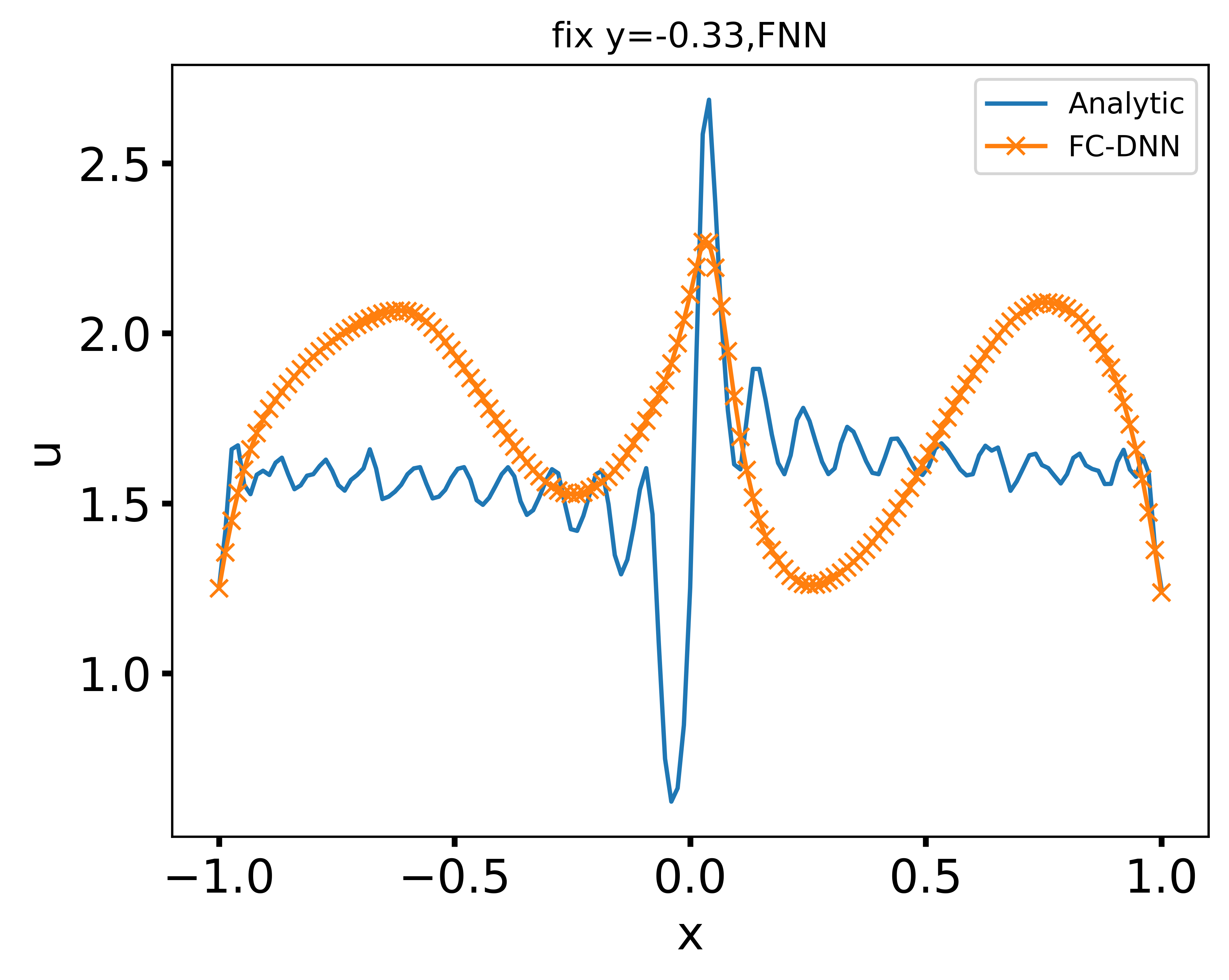}}
\caption{Comparison of a MscaleDNN and a vanilla DNN for an oscillatory solution of the Poisson equation at y=-0.33.
(a) Numerical and exact solution of Poisson equation using MscaleDNN. (b) Failure of a vanilla fully-connected DNN (FC-DNN).}
\label{fig:add_MscaleDNN}
\end{figure*}

To show the effectiveness of MscaleDNN, we solve
 a 2-D Poisson equation with an oscillatory solution in $ \Omega=[-1,1]^2$,
\begin{equation}\label{ex3_Poisson_ritzloss}
        -\Delta u(x, y) =f(x, y),
\end{equation}
where the exact solution is
$  u(x, y) = \frac{1}{N^2} \sum_{m=1}^N \sum_{n=1}^N e^{{\rm sin}(\pi m x)}e^{{\rm cos}(\pi n y)} 
$.
We  choose $N=20$ and use the variational loss \eqref{varloss} for the training.
We compare a MscaleDNN  with a vanilla fully-connected DNN (FC-DNN) for PDE solutions with multiple frequencies.  The following different network structures are used: (1) a FC-DNN with a size 2-3200-3200-3200-1; (2) a MscaleDNN  with scales $\{1,2,4,8,16,32,64,128\}$ for eight subnetworks with a size 2-400-400-400-1 each. The results in Fig.\ref{fig:add_MscaleDNN} clearly show the capability of the MscaleDNN in capturing the high-frequency content of the solution while the FC-DNN failed.

MscaleDNNs have been used for solving various scientific problems. Wang et al. \cite{wang2020multi} apply the MscaleDNN for oscillatory stokes flows in complex
domains. Ying et al. \cite{ying2023multi} use the MscaleDNN-2 structure to construct
multi-scale fusion neural network (MSFN) for solving elliptic interface
problems. Li et al. \cite{li2023solving} utilize MscaleDNN to construct the neural
sparse representation to solve Boltzmann equation with BGK and quadratic
collision model. To solve the stationary nonlinear Navier-Stokes equation, Liu et al. \cite{liu2021linearized} integrate linearizations of the nonlinear convection
term in the Navier-Stokes equation into the training process of MscaleDNN
approximations of highly oscillatory Navier-Stokes solutions.
Chen et al. \cite{chen2025three} applied the MscaleDNN successfully to address the spectral bias in  3-D turbulent wind field reconstruction by leveraging multi-scale PINN, which enhanced the ability to learn from both low and high-frequency data, thus improving the accuracy of flow field predictions.
Also, recently, in \cite{riganti2024effective}, the MscaleDNN has been shown to be critical in addressing numerical error due to the failure of normal neural network in capturing high-frequency multiple-scatttering field by arrays of random dielectric nanocylinders in a study of isotropic SHU (Stealthy Hyperuniform) optical materials, enabling the inverse retrieval of their effective dielectric profiles in a numerical homogenization. The multiscale network architecture proved to be superior over single-scale PINNs, particularly in scenarios with significant multiple scattering effects and high wave numbers. These results highlight the potential applications of the MscaleDNN in real physical systems.

\subsubsection{Spectral bias reduction property of the MscaleDNN, and selection of scales.}

Based on the training dynamics with the NTK framework (\ref{errorEq}), Wang et al. \cite{wang2022convergence} derived a diffusion equation model for the error
evolution in the frequency domain for a learning algorithm by a multiscale neural network in approximating oscillatory functions, assuming a vanishing learning rate and an infinite wide
network in a two-layered network with an activation function $\sigma(x)=\sin x$.

The diffusion model is
\begin{equation}\label{realimagdynamicsbaistrueMsDNN}
	\partial_t \hat{u}^{\pm}(\bm\xi,t)=\nabla_{\bm\xi}\cdot\Big[A_s^{\mp}(\bm\xi)\nabla_{\bm\xi} \hat{u}^{\pm}(\bm\xi,t)\Big]-B_s^{\pm}(\bm\xi)\hat{u}^{\pm}(\bm\xi,t),\quad \bm\xi\in\mathbb R^d,
\end{equation}
where  $\hat u^{\pm}(\bm\xi,\bm\theta(t))$ are the real and imaginary parts of $\hat{u}(\bm\xi,\bm\theta(t))$, respectively, i.e.,
$\hat{u}(\bm\xi,\bm\theta(t))=\hat{u}^+(\bm\xi,\bm\theta(t))+\ri \hat{u}^-(\bm\xi,\bm\theta(t))$,
$u(\bm x,\bm \theta)$ is a zero extension of the error of the network defined by
\begin{equation*}
	u(\bm x,\bm \theta)=\begin{cases}
	0, & \bm x\notin\Omega,\\
	 f(\vx,\vtheta)-f(\bm x), & \bm x\in\Omega,
	\end{cases}
\end{equation*}
$f(\vx,\vtheta)$ is the multi-scale neural network \cite{zhang2023correction} with one hidden layer,
and $A_s^{\pm}(\bm\xi)=\frac{1\pm e^{-2}}{8\pi^2(s+1)}\sum\limits_{j=0}^{s}\alpha_j^{2(d+1)}\widehat{\mathcal G}_j(\bm\xi)$,
$B_s^{\pm}(\bm\xi)=\frac{1\pm e^{-2}}{2(s+1)}\sum\limits_{j=0}^{s}\alpha_j^{2d}\widehat{\mathcal G}_j(\bm\xi)$,
$\mathcal G_p(\bm x)=e^{-4^p|\bm x|^2/2}$, $\bm x\in\mathbb R^d$
is a scaled Gaussian function.

This diffusion model matches well with the error decay of practical training of MscaleDNN in approximating oscillatory functions \cite{wang2022convergence}.

\begin{figure}[ht!]
	\center
	\subfloat{\includegraphics[scale=0.21]{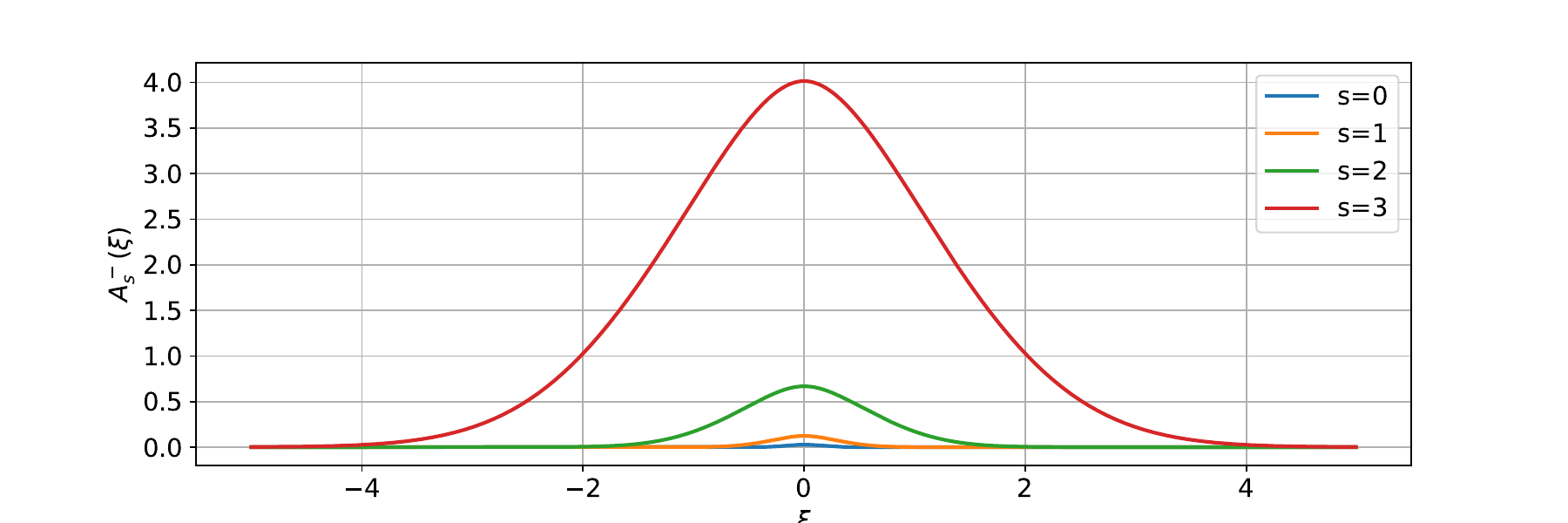}}
	\subfloat{\includegraphics[scale=0.21]{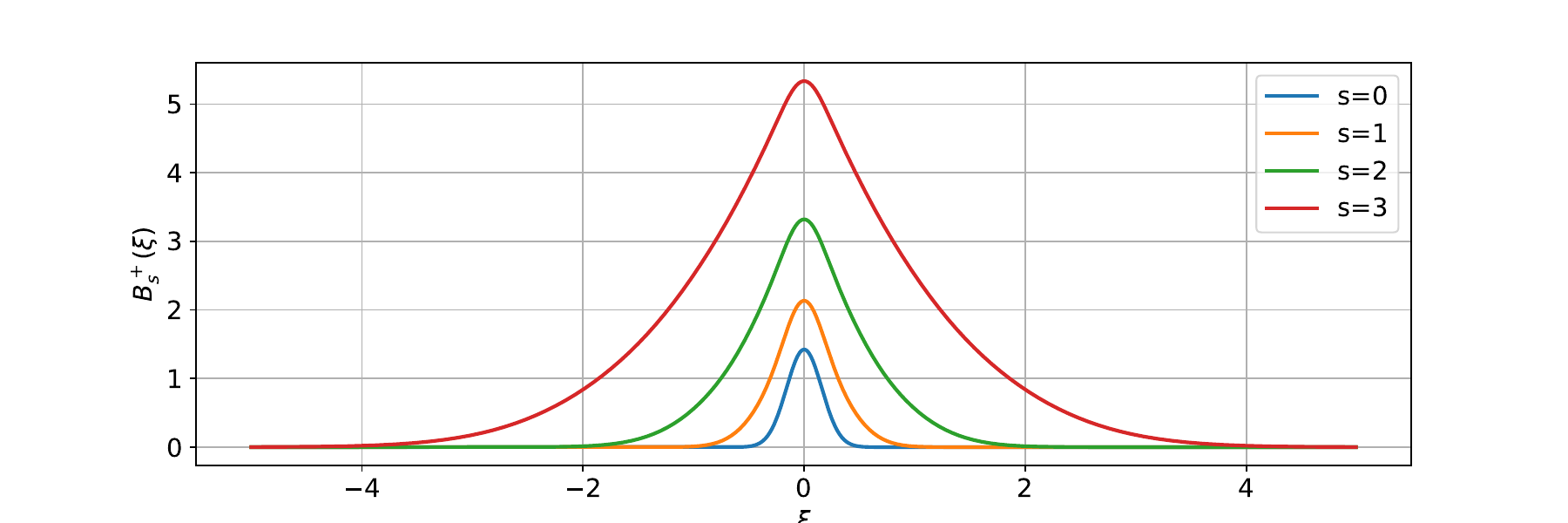}}\\
	\subfloat{\includegraphics[scale=0.21]{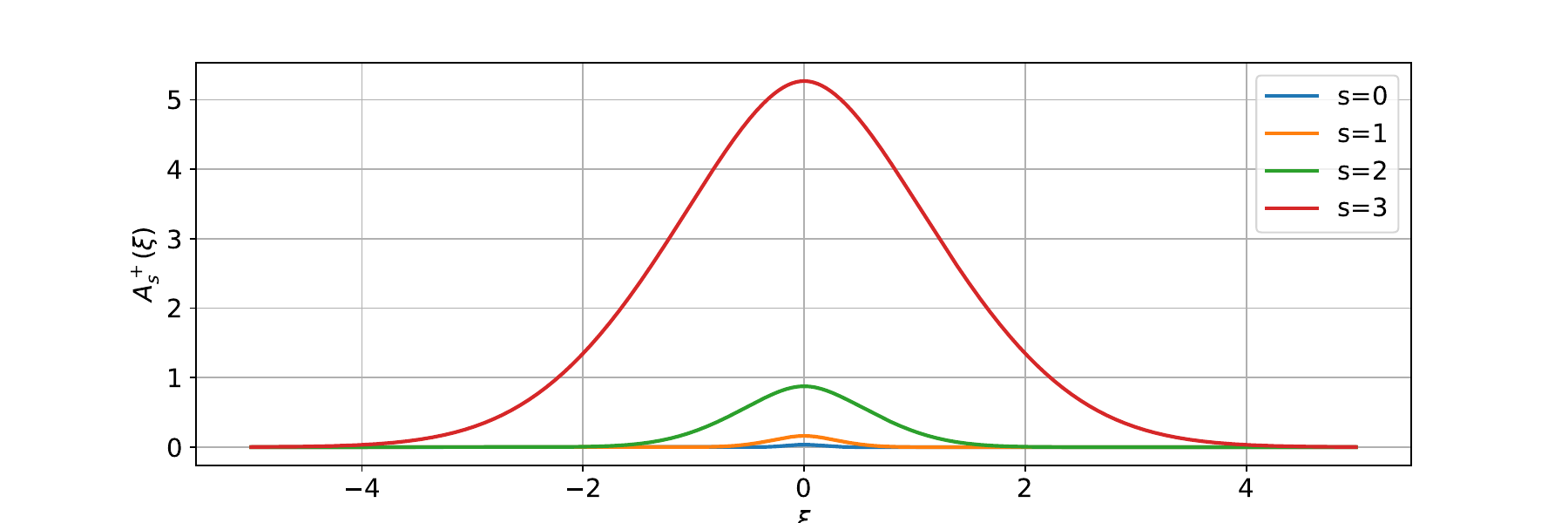}}
	\subfloat{\includegraphics[scale=0.21]{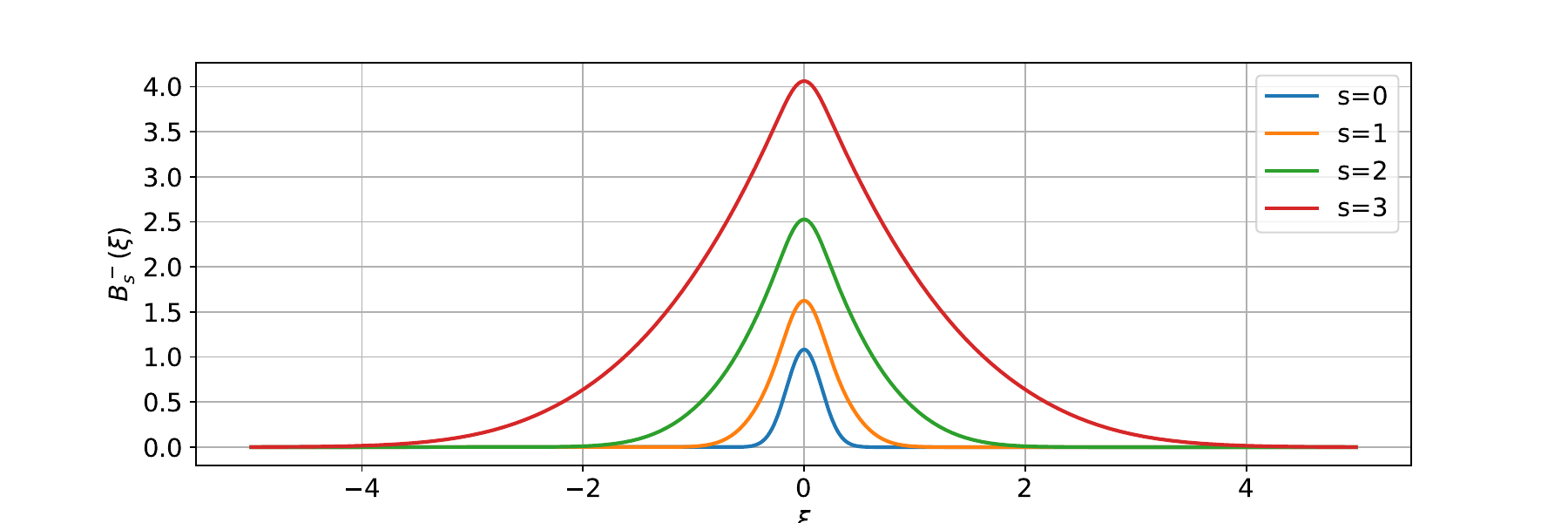}}
	\caption{Diffusion coefficients $A^{\mp}_s(\xi)$ (left) and $B^{\pm}_s(\xi)$ (right) with $\alpha_j=2^j$,  $s=0, 1, 2, 3$.}%
	\label{coefficientsplot}%
\end{figure}


A sample of $\{A_s^{\mp}(\xi), B^{\pm}_s(\xi)\}$ for a 1-dimensional case are plotted in Fig. \ref{coefficientsplot} for up to 4 scales two-layer MscaleDNNs. We can see that both $A_s^{-}(\xi)$ and $B^{+}_s(\xi)$ have larger support in frequencies with  increasing of $s$, thus providing the mathematical explanation for the capability of the MscaleDNNs in reducing spectral bias as shown by extensive numerical tests. Also for $s=0$, this result confirms the spectral bias of the two-layered fully-connected neural network with a sine-activation function.

The number and specific values of the scales used in the MscaleDNN should depend on the target solution to be learned, therefore, in principle, they could be made as hyper-parameters of the MscaleDNN. Indeed, in a recent work \cite{zhouFreqadpt2024}, a frequency (i.e., scale) adaptive PINN was investigated, which shows that dynamically selected scales could improve much the DNN solution for PDEs with oscillatory solutions.

\subsection{Fourier feature networks}



The first hidden layer in a neural network can be viewed as a mapping that lifts the input to another space, often of a higher dimension than that of the input space. In the MscaleDNN framework, the first layer lifts the input in $d$ dimensions to $m_1$ dimensions. Another lifting operation can be done by setting all trainable weights in the first hidden layer to fixed values. For example, as in \cite{tancik2020fourier}, the Fourier Feature Networks set all the weights to fixed values and biases to zero with a sinusoidal activation function, i.e.,
\begin{equation}
   \gamma(v) =[ a_1 cos(2\pi b^T_1 v), a_1 sin(2\pi b^T_1 v), \cdots, a_m cos(2\pi b^T_m v), a_m sin(2\pi b^T_m v)]^T,
   \label{Ffeature}
\end{equation}
where  $v$ is the input, and $a_k$'s and $b_k$'s are fixed non-trainable parameters.  In particular, in the positional encoding of Transformer \cite{vaswani2017attention} $a_k=1$ for all $k$ and $b_k = c^{k/m}$ with a scale $c$,  similar to the scales $\alpha_i$ in the MscaleDNN (\ref{f_app}). In other cases, $b_k$  can be sampled from a Gaussian distribution. In \cite{mildenhall2021nerf}, a multiscale input structure $\gamma(v)$  \eqref{Ffeature}    was used with $a_k=1, b_k=2^k$  in representing scenes as neural radiance fields for view synthesis, achieving a great success in computer vision.
The resulting network is exactly the MscaleDNN with activation function $\sigma(x)=sin x$, thus the analysis in Section 4.2.1 also applies to the Fourier feature network.

An extension of  the Fourier feature networks was given in \cite{wang2020eigenvector} to solve PDEs by setting the coefficients in (\ref{Ffeature}) $a_k=1$ and sampling each $b_k$ from Gaussian distribution $N(0,\sigma_k^2)$ with variance $\sigma_k$.
Another transform was used in multi Fourier feature networks (MFFNN) \cite{li2023deep}
by setting $\vgamma(\vx) =[ cos(2\pi \vB \vx), \vx, sin(2\pi \vB \vx)]^T$ with entries of $\vB \in R^{m\times d}$ (m is a given integer) sampled from an isotropic Gaussian distribution $N(0,\sigma^2)$ and the standard deviation $\sigma$ treated as a hyper-parameter. Combination of different choices of $\gamma(v)$ has been used in \cite{guan2023mhdnet} by using $\gamma(\vx) = [cos(2\pi \vB \vx), \alpha * \vx, sin(2\pi \vB \vx)]^T$ for the simulation of MHD problems.

\subsection{Activation function optimization}

The activation functions play an essential part in the learning process of artificial neural networks. There are many choices of activation functions in deep learning, such as ${\rm ReLU}(x)$, ${\rm tanh}(x)$, ${\rm GeLU}(x)$ \cite{hendrycks2016gaussian},  rational activation functions \cite{boulle2020rational}, and sigmoid type activation functions etc.
Empirical works also find that using B-spline \cite{wang2024expressiveness}
or sinusoidal activation function with proper initialization \cite{sitzmann2020implicit} can also mitigate the spectral bias. A survey of activation functions in multi-scale neural networks is referred to Jagtap et al. \cite{jagtap2023activationsurvey}.
To improve the effectiveness and accuracy of deep neural networks (DNNs) in approximating functions with multiple frequencies and to mitigate the impact of spectral bias, one of the strategies involves optimizing the choice of activation functions, pursued in the following works.

\subsubsection{Adaptive activation function}

Jagtap et al. \cite{jagtap2020adaptive} introduced an adaptive activation function in the form $\sigma (n a \vx)$ with an integer $n$ and a trainable hyper-parameter $a$ in the NN. With a large $n$, this effectively increases the sampling range of the weights and the bias in the network structure while the trainable $a$ introduced another degree of freedom for the training of the network parameters. The resulting adaptive activation function NN will take the following form
\begin{equation}
    f_{\vtheta}(\vx) = \vW^{[L]} \sigma\circ (na(\cdots (\mW^{[1]} \sigma\circ(na(\mW^{[0]} (\vx) + \vb^{[0]}) ) + \vb^{[1]} )\cdots))+\vb^{[L]}. \label{adaptive activation function}
\end{equation}

To some extent, this factor $na$ plays a similar role of the scales in the MscaleDNN, so some reduction of spectral bias is expected and numerical results in \cite{jagtap2020adaptive} validate this approach for both supervised learning and PDE solving up to a certain high frequency, such as $\sin (12x)$.



\subsubsection{Reproducing  activation functions}
Another way to improve the performance of the DNN is to enrich the activation function choices using a collection of activation function form so that traditional objects in approximation theory, such as polynomials, Fourier basis, wavelet basis can be reproduced by using enriched activation functions.
Liang et al. \cite{liang2021reproducing} proposed such an approach with reproducing activation functions (RAFs), which employ several basic activation functions and their learnable linear combinations to construct activation functions for different neurons.
For the $i$-th neuron of the $l$-th layer, the formula of the neuron-wise RAF is given by
\begin{equation}
    \sigma_{i,j}(x)=\sum_{p=1}^{P}\alpha_{p,i,l}\gamma_p(\beta_{p,i,l}x),
\end{equation}
where $A = \{\gamma_1 (x), \cdots, \gamma_P (x)\}$ is a set of $P$ basic activation functions,
$\alpha_{p,i,l}$  and $\beta_{p,i,l}$ are learnable  combination coefficient and scaling parameter, respectively.
Numerical results have shown better performance with smaller number of parameters for audio and image restoration, PDEs and eigenvalue problems, such as solving a PDE with a frequency of $\sin(2\pi x_1)\sin (2\pi x_2)$, compared with traditional NNs.
In particular, for solutions with both low and high frequencies, Liang et al. \cite{liang2021reproducing}  used the poly-sine-Gaussian network with
$ A = \{x, x^2, sin(\beta x), e^{-(\beta x)^2}\}$,  in which $\beta$ is  a trainable hyper-parameter.

\subsubsection{Multiple activation neural network}
Compared with the enrichment of activation function at the neuron level by Liang et al. \cite{liang2021reproducing},
Jagtap et al. \cite{jagtap2022deep} proposed an enrichment of the network using a combination of a collection of activation functions $\{ \phi_k(\vx)\}$ for each neuron. A two layered of the so-called Kronecker neural networks (KNNs)  takes the following form
\begin{equation}
    f(\vx, \theta) = \sum_{i=1}^m \left [ \sum_{k=1}^K \alpha_k \phi_k(\omega_k(w_ix+b_i))  \right],
\end{equation}
where each of the activation function $\phi_k(\vx)$ shares the same weights and bias and the network parameters are $\theta= \{ c_i, w_i, b_i \}_{i=1}^m \cup\{\alpha_k, \omega_k \}_{k=1}^K$.

A particular choice of the activation function is to set $\phi_1$ to be any standard activation function such as ReLU, tanh, ELU, sine, Swish, Softplus etc., and the remaining activation functions are $\phi_k(\vx)=n\sin((k-1)nx)$ or $n\cos((k-1)nx)$, $2\leq k \leq K$ resulting in a Rowdy-Net  \cite{jagtap2022deep}. This choice is similar to the MscaleDNN of \eqref{f_app} in the first hidden layer with a sinusoidal activation function $\sigma(x)=sinx$ and scaling factor $\alpha_k=k$.

Jagtap et al. \cite{jagtap2022deep} also provided a theoretical analysis which reveals that  under suitable conditions,  the Kronecker network induces a faster decay of the loss than that by the feed-forward network at the beginning of the training.



\subsection{Random Feature Method}


Chen et al. \cite{chen2022bridging} proposed a random feature method (RFM) for solving PDEs, a natural bridge between traditional and machine learning-based algorithms.
It  adds several additional components including multi-scale representation and rescaling the weights in the loss function.


Following the random feature model in machine learning,  one can construct an approximate solution $u_M$ of $u$ by a linear combination of $M$ network basis functions $\{\phi_m\}$ over $\Omega$ as follows
\begin{equation}
    u_M(\vx) = \sum_{m=1}^M u_m \phi_m(\vx).
\end{equation}
Generally speaking, the basis functions will be chosen as the ones that occur naturally in neural networks, for example:
\begin{equation}
    \phi_m(\vx) = \sigma(\vk_m \cdot \vx + b_m),
\end{equation}
where $\sigma$ is some scalar nonlinear function, $\vk_m$, $b_m$
are some randomly selected but fixed parameters.

\paragraph{Partition of unity and local random feature models}

Random feature functions are globally defined, while the solution of a PDE typically has local variations, possibly at small scales.
To accommodate this, RFM constructs many local solutions, each of which corresponds to a random feature model, and piece them together using partition of unity (PoU). The PoU starts with a set of points $\{\vx_n \}^{M_p}_{n=1} \subset \Omega$, each of which serves as the center for a component in the partition. For each $n$, a normalized coordinate is defined
\begin{equation}
    \Tilde{\vx}_n = \frac{1}{\vr_n}(\vx-\vx_n), \quad n = 1,\cdots,M_p, \label{eq:rfm_rn}
\end{equation}
where $\{ \vr_n \}$ is pre-selected, and, similar to the scaling in MscaleDNN, the coefficient $r_n$ in Eq. (\ref{eq:rfm_rn}) plays an important role in learning a wider frequency range.
Next, for each n,  $J_n$ random feature functions are constructed
\begin{equation}
   \phi_{n,j}(\vx) = \sigma(\vk_{n,j} \cdot \Tilde{\vx}_n + b_{n,j}), \quad j=1,\cdots,J_n,
\end{equation}
where the feature vectors $\{(\vk_{nj},b_{nj})\}$ often chosen randomly and then fixed.
In this way the locally space-dependent information is incorporated into $M = \sum^{M_p}_{n=1} J_n$  random feature functions.
 As for the construction of the PoU,
 one-dimensional PoU can be defined as  \begin{equation}
    \psi_n(x)=\mathbb{I}_{-1\leq
\Tilde{x}_n}
    \leq 1
 \end{equation}
  or
  \begin{align*}
 \begin{split}
    \psi_n(x)
     = \left \{
 \begin{array}{ll}
     \frac{1+sin(2\pi\Tilde{x})}{2},                    &
     -\frac{5}{4}\leq
     \Tilde{x}_n
<-\frac{3}{4}
     \\
     1,     &
     -\frac{3}{4}\leq
     \Tilde{x}_n
<\frac{3}{4}
     \\
     \frac{1-sin(2\pi\Tilde{x})}{2},                    & \frac{3}{4}\leq
     \Tilde{x}_n
     \\
     0, & otherwise.
 \end{array}
 \right.
 \end{split}
 \end{align*}
where the relation between $x$ and $\Tilde{x}_n$ is given by Eq. (\ref{eq:rfm_rn}).
High-dimensional PoU can be constructed using the tensor product of one-dimensional
 PoU functions, i.e., $\psi_n(\vx) = \prod_{k=1}^d \psi_n(x_k)$, where $\vx=(x_1,x_2,\cdots,x_d)$ and $x_k \in \mathbb{R}$ for each $k=1,2,\cdots,d$.

The RFM approximate solution $u_M$ is given by
\begin{equation}
    u_M(\vx)=\sum_{n=1}^{M_p} \psi_n(\vx) \sum_{j=1}^{J_n} u_{nj}\phi_{nj}(\vx),
    \label{random_feature_2.8}
\end{equation}
which resulted in a over-determined system  to be solved by some QR-based least square solvers.

\paragraph{Multi-scale basis}
In some situations, (\ref{random_feature_2.8}) alone is less efficient in capturing the large scale features in the solution. Therefore, on top of the PoU-based local basis functions, one can add another global component:
\begin{equation}
    u_M(\vx) = \sum_{m=1}^M u_m \phi_m(\vx) + \sum_{n=1}^{M_p} \psi_n(\vx) \sum_{j=1}^{J_n} u_{nj}\phi_{nj}(\vx).
\end{equation}

This RFM  is shown to be able to solve a low-dimensional PDE problem with wide frequency range with an almost machine accuracy even in a complex geometry region \cite{chen2022bridging}. Since the RFM utilizes the random feature basis, it is a mesh-less method, therefore, it can handles problems with complex regions. The machine accuracy can be achieved because the problem is not solved by gradient based methods but solved by a matrix inversion method with QR decompositions.

\section{Hybrid approach}

Classic iterative solvers excel at reducing high frequency errors while the normal DNNs are effective in lower frequency ones, therefore, it is natural to combine them to arrive at a solver with uniformly fast error reductions across all frequencies, giving rise to hybrid methods.

\subsection{Jacobi-DNN}
A natural idea of frequency decomposition is to utilize a hybrid approach of using conventional iterative method such as Jacobi or Gauss-Seidel methods for the high-frequency error reduction and a DNN-based method for the fast convergence of lower frequency.






Xu et al. \cite{xu2019frequency} combined DNNs with
conventional iterative methods to accelerate the convergence of both low and high
frequencies for computational problems. Take Poisson's
equation as an example. First, the DNN was used to solve the Poisson's
equation with a certain optimization
steps (or epochs). Then,  the Jacobi
method with the initial value obtained from the trained DNN is carried out to eliminate the error in high frequencies using the smoothing properties of the iterative method.

Huang et al. \cite{huang2020int} developed a similar approach (Int-Deep), applying a deep-learning-initialized iterative method and demonstrated that the Int-Deep method is effective for low-dimensional PDEs, such as semilinear PDEs, linear and nonlinear eigenvalue problems etc. Chen et al. \cite{chen2023properties} used a Broad Learning System (BLS), which learns low frequency faster similar to DNNs, and found  that low and high frequencies can be learned fast for low-dimensional problems by a BLS-Jacobi method.

However, since such hybrid methods use conventional methods in the second stage, it also suffers similar difficulties as conventional methods. For example, such hybrid methods can only be applied to low-dimensional problems.
\subsection{HINTS}

 Zhang et al. \cite{zhangGK2024} proposed a hybrid approach of integrating neural operator DeepONet and standard relaxation methods, yielding a hybrid iterative numerical transferable solver (HINTS).
Firstly, a DeepONet  needs to be trained  offline before employing HINTS to approximate the solution operator of a linear differential equation.
Then, HINTS starts by discretizing the linear differential equation and  alternately adopts the
numerical iterator and the DeepONet iterator with ratio $1 :(n_r-1)$ (i.e., DeepONet with
proportion $1/n_r$), until the iterative solution converges.
Instead of existing traditional solvers which adopt a fixed relaxation method in each iteration,  HINTS creates a second branch of the iterator using a
trained DeepONet.  The numerical iterator may be chosen from  traditional solvers, such as the Jacobi method, the Gauss-Seidel method et., which  are known to be efficient for high-frequency modes but not for low-frequency modes;
the DeepONet solver provides fast approximate solutions for low frequencies but may contain polluted solutions at the high frequencies. A proper combination of these two methods enhances the convergence and enables fast and uniform convergence across all frequencies.
 Such a combination of solvers incorporates
advantages from both worlds, improving the performance beyond individual solvers: classical solvers either are slow or even fail to convergence; an individual deep learning solver
may only provide an approximate solution.

HINTS exploits such a bias to tackle low-frequency modes that are otherwise difficult for classical solvers which, essentially, are biased towards high
frequencies. By replacing the classical solver with the DeepONet solver for
a limited proportion of iterations, HINTS balances the convergence rates
across the spectrum of eigenmodes, significantly alleviating the spectral bias
and improving computational efficiency.

\subsection{MgNet, Multi-level, Multi-grade Nets}

When discretized, a linear PDE is transformed into a linear system of equations, represented as
\begin{equation}
    A_{\eta}u=f. \label{eq:au=f}
\end{equation}
In the context of solving large-scale linear systems, iterative methods are frequently employed as an alternative to direct matrix inversion. As noted before, Jacobi and Gauss-Seidel iterations exhibit slow convergence for low-frequency components. A  classical iterative approach for rapid convergence across all frequency regimes is the multi-grid (MG) method \cite{hackbusch2013multi}. The MG method comprises two distinct stages: smoothing and coarse grid correction. During the smoothing stage, high-frequency errors are eliminated through the application of smoothers, such as Jacobi or Gauss-Seidel iterations.

Crucially, MG methods leverage a coarse grid solution to eliminate low-frequency errors. As demonstrated in a prior study by He et al. \cite{he2019mgnet}, the operations involving prolongations, restrictions, and select smoothers within the MG method framework can all be formulated as convolution operations. Subsequently, through the substitution of these operations with convolutional neural networks with trainable parameters, the Multigrid Network (MgNet) framework for solving PDEs has been developed \cite{chen2022meta}. MgNet is a DNN-based iterative approach, and its loss function is derived from the least squared error of Eq. \eqref{eq:au=f}, with a single MG iteration step to approximate the solution $u$.

In MgNet, the solution domain is dissected into grids with both fine and coarse resolutions. The neural network embedded within this framework serves to establish a connection between solutions respectively spanning low and high-frequency ranges. Consequently, MgNet effectively harnesses the strengths of the traditional MG method of frequency uniform convergence,  alongside the approximation capabilities intrinsic to neural networks.
However, due to the use of grids, the MgNet will be limited to solving PDEs in low dimensions.

Another multi-level learning by neural network is an hierarchical approach proposed in \cite{han2021hierarchical}, which reduces the residual of the numerical solution using a series of networks recursively. The first neural network is used to learn the target function with a given error loss, such as the least squared error. Once the training stalls, the parameters of the first network are fixed. Then, a second neural network, which takes the same input as the first network, is optimized to reduce the residual error between the target function and the sum of two networks. Recursively, the k-th network, which takes the same input as previous (k-1) networks, is optimized to reduce the residual error between the target function and the sum of all k networks. Therefore, each additional network acts as an incremental correction to previous networks, similar to the Mg-Net approach. To further improve accuracy, \cite{aldirany2024multi} assigns appropriate coefficients for each network.

Multi-grade learning model is an alternative approach to realize a multi-stage learning for different frequencies  \cite{xu2023multi}, where the k-th network takes the output of the (k-1)-th network as input instead of the original input, and the output of the k-th network is always trained to fit the target output. The training of the k-th network is based on the already trained (k-1)-th network. Such multi-grade learning has been empirically shown to be faster than one-grade learning for multi-frequency problems.

\section{Multi-scale neural operators and diffusion models}
As discussed previously, an alternative approach to solve PDEs is to learn the PDE operator, which represents the relation between the solution and some physical quantities such as the material properties (coefficients in the PDEs), forcing, or initial or boundary conditions.  For instance, for wave scattering, operator learning could be used to establish the relation between the scattering field and the conductivity or permittivity of the scatterer or the incident waves or both. Once such an operator is learned, as an immediate application, it could be used to quickly predict the scattering field for other un-seen scatterers or incident waves. More importantly, the learned operator could be used to solve the challenging problem of inverse medium problems for imaging applications in medical and material sciences and geophyiscs.
On the other hand, if a neural network is used directly to parameterize the solution for the problem, we will need to re-train the neural network for different materials or boundary conditions.  In training a network to learn the PDE operator, the network takes a function as input and outputs a function. For the input function, it could be represented by its values on fixed grids. For the output function, there are two common ways, one is to output the function values on a group of fixed grids, the other is that the network also takes a position $x$ as input in addition to the input function and outputs the target function at $x$. Such neural network is commonly referred to as neural operator, such as DeepONet \cite{lu2021learning}, FNO \cite{li2020fourier} etc.
Neural operators such as MgNet and UNet-based DNNs  \cite{gupta2022towards,rahman2022u,ovadia2023ditto}
 use a multi-scale sampling to efficiently solve functions with multiple frequencies.

It has also been observed that a neural network representation of a PDE operator could also suffer from a similar spectral bias for high frequency functions. For example, consider a functional mapping from $a(x)$ to $u(x)$ in the following simple form,
\begin{equation}
    u(x)= \mathcal{G}[a(x)] (x):= \sin (m a(x)).
\end{equation}
Even the input function $a(x)=a$ is a constant, for large $m$, the mapping $\mathcal{G}$ will have a high frequency dependence on $a$ as the latter varies.

To alleviate the spectral bias, multiscale neural operators have been proposed by using scaled input to a normal neural operator such as multiscale FNO (MscaleFNO) \cite{you2024} and multiscale DeepOnet  \cite{LiuCai2021}.


For instance, in \cite{you2024}, the 1-D Helmholtz equation for the scattering field with a Dirichlet boundary condition was considered
\begin{equation}
    \begin{cases}
    u'' + (\lambda^2 + c\omega(x))u = f(x) & x \in [-L,L], \\
    u(-L) = u(L) = 0,
\end{cases}
\label{helmholtz}
\end{equation}
where $\lambda$ is the wave number of the homogeneous background material and $c \omega(x)$ corresponds to the variable wave number perturbation of the scatterer, $u(x)$ is the scattering field and $f(x)$ is the source term. We are interested in learning the operator
\begin{equation}
    \mathcal{G}: \omega(x) \mapsto u(x),
\end{equation}
for a fixed source term $f(x)$ in the form of
\begin{equation}
        f(x) = \sum_{k=0}^{10} (\lambda^2-\mu_k^2) \sin(\mu_k x), \; \mu_k = 300 + 35k.
        \label{e11}
\end{equation}
A MscaleFNO of the following form was used \cite{you2024} to predict the scattering field from $\omega(x)$ for a fixed source $f(x)$ for the case $L=11,\lambda = 2, c = 0.9\lambda ^2 = 3.6$,
\begin{equation}
u\left(x\right) = \sum_{i=1}^{N} \gamma_i \operatorname{FNO}_{\theta_m}\big [c_i x, c_i\omega(x)\big](x).
\label{msFNO}
\end{equation}

\begin{figure}[h]
\centering
\subfloat[Normal FNO]{\includegraphics[width=0.45\linewidth]{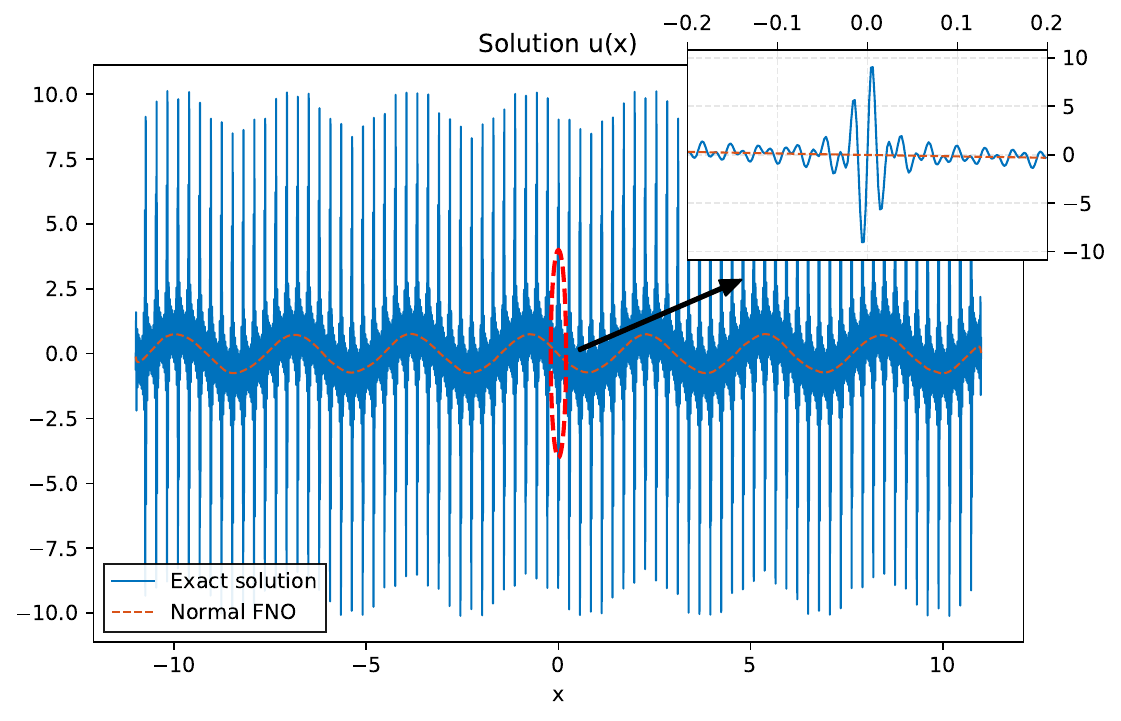}}
\subfloat[MscaleFNO]{\includegraphics[width=0.45\linewidth]{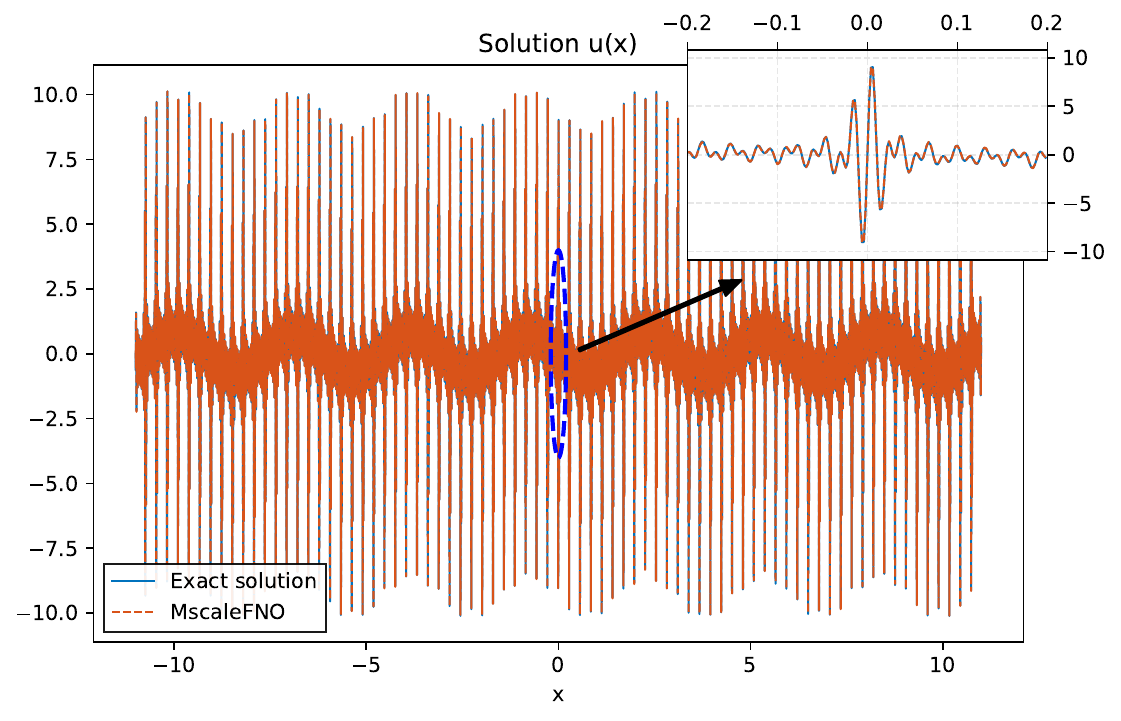}}
\caption{Performance comparison of a Normal FNO and a MscaleFNO in predicting highly oscillatory solution $u(x)$ of a 1-D Helmhotlz equation from its coefficient $\omega(x)$  (a) Prediction results of the normal FNO compared with the exact solution, with an inset showing details in the region $[-0.2,0.2]$ and the failure of FNO to capture the oscillatory behavior of the solution; (b) Prediction results of the MscaleFNO, showing its capability in capturing both global behavior and local details of the solution, as depicted in the zoomed-in insert.}
\label{helm_io}
\end{figure}

\begin{figure}[h]
\centering
\includegraphics[height=4 cm,width=6cm]{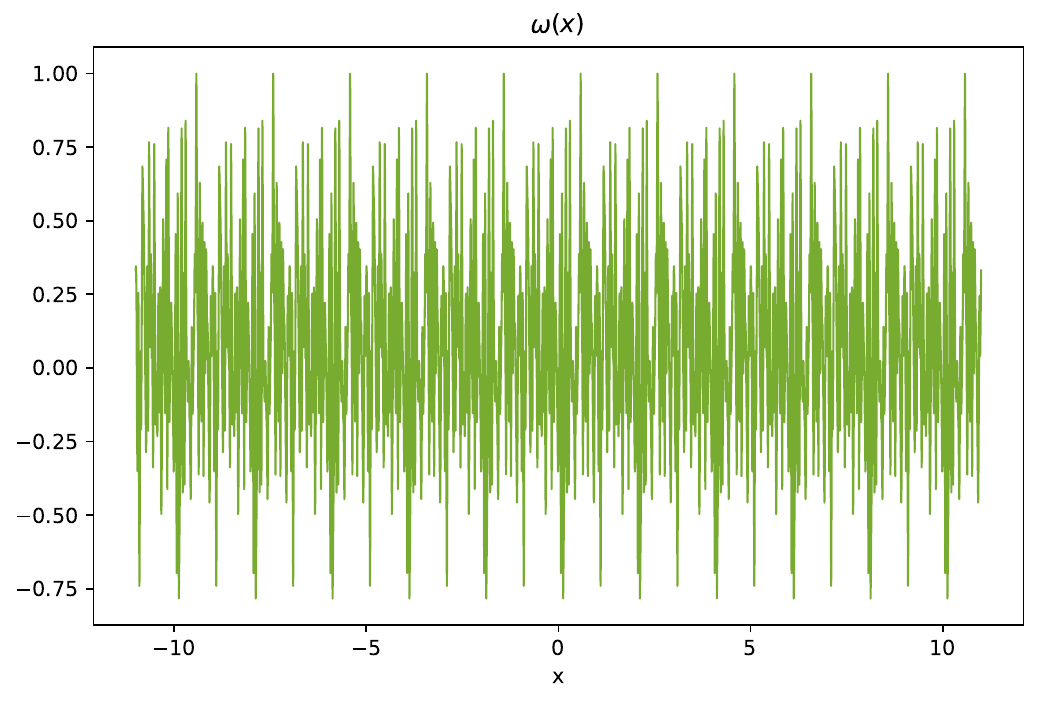}
\includegraphics[height=4 cm,width=6cm]{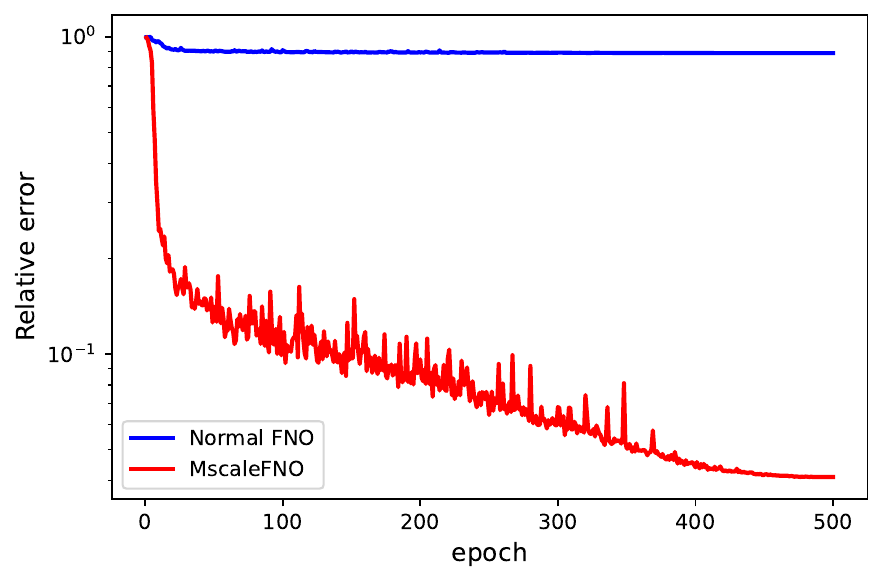}
\caption{(Left) The input function $\omega(x)$ with rapid oscillations across the domain; (Right) Convergence comparison of relative errors during 500 training epochs between Normal FNO and MscaleFNO. The MscaleFNO demonstrates superior learning capability with the relative error decreasing to approximately $4\times10^{-2}$, while Normal FNO remains at a higher error level around $8\times10^{-1}$.}
\label{helm_ec}
\end{figure}

Fig. \ref{helm_io} shows the predicted scattering field by a normal FNO (left panel, only an average value of the solution is predicted) and a multiscale FNO (right panel) with
$N=8$ scales with $c_i = 1, 80, 160, 200, 240, 280, 360$, and $400$ and similar total number of parameters as the normal FNO. Fig. \ref{helm_ec} shows the profile of the input function $\omega(x)$ (left panel) and the comparison of the convergence history of the FNO and the MscaleFNO after 500 epochs of training (right panel).

In a recent work \cite{HANO2024}, an hierarchical attention neural operator was proposed to address the spectral bias issue by leveraging diffusion models in generative AI.
Diffusion models \cite{ho2020denoising} and GANs \cite{goodfellow2020generative} were developed for various tasks such as image and video generation.
In \cite{oommen2024}, a diffusion model is integrated with various neural operators to improve their spectral representation of turbulent flows.
In \cite{molinaro2024generative, oommen2024, wu2024highflexibility}, the authors leveraged generative modeling techniques to learn PDE operator with a wide range of frequencies which can also  mitigate the spectral bias typically exhibited by the DNNs.

\section{Conclusions and open problems}
\label{sec:conclusions}

We have provided a survey of the current efforts in developing various methods to overcome the spectral bias of deep neural network in learning PDE solutions with wide frequency contents. This line of research is generating many new results at a fast pace, this survey will surely miss some worth developments. Most of the methods reviewed in this survey in fact can trace their roots in classical approximation theories and traditional numerical methods, such as the multi-resolution idea of wavelets (multiscale DNN, Fourier features, adaptive activation functions), wave-ray MGM \cite{livshitsRayMG2006} (PhaseDNN), domain decomposition (Random feature method), multigrid (Mg-Net), multi-level approaches. The challenge is how to develop effective mathematical tools  suitable for the compositional structure of the DNNs  to investigate the effects of various constructs of DNN on removing or mitigating the effects of spectral bias, so fast learning is not limited to only lower frequency component of the solution while stalling or even failing at higher frequencies. Thus, it is of great relevance to understanding the dynamics of the network weights during the training and their relation to the interpolative and expressivity power of the network. The recent results on information bottleneck and condensation of network weights are a good starting point in gaining insights of the mechanism of DNN learning, thus providing researchers new inspirations of developing more sophisticated DNN models as well as analysis tools to  best resolve the spectral bias of DNN and to achieve a robust and mature way using DNN as a PDE solver, especially for high dimensional PDE problems or complex geometry problems.


\textbf{DNN models and training process} Over-parameterized neural networks have numerous global minima, which may have different generalization performance. The characteristic of the training process of neural network is critical to understand what kind of minima the training will find. For a normal neural network, the spectral bias leads to a minima of low-frequency interpolating capability. There are various phenomena found in the training process, which may give indication if such a spectral bias may occur.

In \cite{tishby2015deep}, a phase transition phenomena has been observed in light of information bottleneck theory, which shows two phases of training a neural network, a fitting phase and a diffusion one. In the fitting phase, the training is primarily driven to decrease the loss and the mutual information \cite{duncan1970calculation} between the input and output significantly increases. In the diffusion phase, however, noise takes a dominant role in the training and weights evolve chaotically, resembling a diffusion process. This phenomenon is evident in the training of plain PINN \cite{anagnostopoulos2023residual} and several synthetic examples were also employed in \cite{xu_training_2018}  to demonstrate that the fitting process corresponds to the learning of low-frequency content while the diffusion one to the high-frequency content. A future line of work may be to examine various DNN models reviewed in this survey through the information bottleneck theory to see how their phase transitions are possibly related to their effectiveness in overcoming spectral bias.

Another informative phenomenon during the training of the network is the condensation of neurons during the training process, i.e., the weights of hidden neurons in two-layer $\rm{ReLU}$ neural networks are shown to condense into finite orientations during training in the non-linear regime, particularly with small initialization \cite{luo2021phase,zhou2022empirical} or with dropout regularization \cite{zhang2022implicit,zhang2023stochastic}. In the early training stage, neurons in the same layer would condense onto very few directions \cite{zhou2022towards,zhou2023understanding,chen2023phase}. Such condensation onto few directions plays a critical role that resets the neural network to a simple state, whose expressivity is small to fit only very low-frequency content. Then, due to non-zero loss, the network would condense on more and more directions in order to increase the expressivity. Such condensation with more and more directions can ensure the network to fit the target data points but with as low complexity as possible. How new DNN models would affect the condensation process during the training is also an interesting topic for understanding how they alleviate the spectral bias or designing more efficient structure to accelerate the training.

\textbf{Frequency domain convergence behaviors of various DNN models.} The convergence and generalization error analysis of various DNNs for solving PDEs have been investigated extensively recently on PINNs
\cite{shinn2020,shinn2023,mishra2022a}, DeepRitz type variational DNNs \cite{lu2022,lu2021,dondl_DRM2022}, deep Galerkin DNNs \cite{yang_DGM2023}, and DeepOnet \cite{shinn2021,mishra2022b}, just to list a few. However, not much work is done to illuminate the frequency dependence of the convergences in terms of spectral bias free behaviors. The loss functions for training DNN come from different mathematical formulations of the solution approach for the PDEs, including the residual based PINN, variational energy based DeepRitz, Garkerkin methods, stochastic differential equation based methods, etc. The specific form of the loss function will require a specialized tool to analyze the convergence and exam how the spectral bias manifests itself.
The convergence analysis for fitting problems from the
perspective of frequency domain has been rigorously analyzed in the linear regime
of normal fully-connected neural networks  \cite{zhang2021linear,luo2022exact,basri2019convergence,cao2019towards} or a two layer multiscale networks \cite{wang2022convergence}. Such analysis for PINNs in the linear regime can be anticipated in future works, especially for solving PDEs.
Extending Fourier analysis from linear regime to nonlinear regime has intrinsic difficulty due to the composition of functions. Some qualitative analysis can be found in \cite{luo2019theory}. However, quantitative analysis remains unsolved. Moreover, finding the appropriate tools for the convergence study for DNN using other form of loss functions and training strategies such as the min-max optimization, SDE based DNNs {\color{blue} as well as PDE operator learning}, remains open problems.


Other challenges in making DNN based algorithms to be free of spectral bias and achieve robust and frequency uniform convergence as a practical simulation tool exist, including how to select various hyper-parameters, including penalty constants in various terms of the loss functions, activation function, and learning rate, etc. Those issues are not addressed in this review. However, with recent progress in reducing the adverse effect of the spectral bias of DNN in  approximating PDEs solutions, the opportunities for research in the DNN based computational methods for solving PDEs are numerous and rewarding, and the fundamental progress on algorithms and mathematical analysis  will have far-reaching impact in many application areas.

\section*{Acknowledgement}
Z.X. acknowledges the financial support provided by the National Key R\&D Program of China  Grant No. 2022YFA1008200, the National Natural Science Foundation of China Grant No. 92270001, 12371511, Shanghai Municipal of Science and Technology Major Project No. 2021SHZDZX0102, and the HPC of School of Mathematical Sciences and the Student Innovation Center, and the Siyuan-1 cluster supported by the Center for High Performance Computing at Shanghai Jiao Tong University. W.C. acknowledges the financial support provided by the US National Science Foundation grant DMS-2207449.



\bibliographystyle{elsarticle-num}


\end{document}